\documentclass{article}
\usepackage{amssymb,amsmath,vmargin}
\title{
\vspace{.5ex}
\begin{center}
\LARGE {\bf Sojourn time in an union of intervals for diffusions}
\end{center}\vspace{.21ex}}

\author{Aim\'e LACHAL
\footnote{
\mbox{Postal adress: \textsc{Institut National des Sciences Appliqu\'ees de Lyon}}
\hspace*{10\textwidth} \mbox{}
\mbox{P\^ole de Math\'ematiques/Institut Camille Jordan}
\hspace*{10\textwidth} \mbox{} \mbox{B\^atiment L\'eonard de Vinci, 20 avenue Albert Einstein}
\hspace*{10\textwidth} \mbox{} \mbox{69621 Villeurbanne Cedex, \textsc{France}}
\hspace*{10\textwidth} \mbox{} E-mail: {\tt aime.lachal@insa-lyon.fr}
\hspace*{10\textwidth} \mbox{} Web page:
{\tt http://maths.insa-lyon.fr/$\mbox{}^{\sim}$lachal}
}
\\[2ex]
\begin{small}
\textsl{Universit\'e de Lyon, CNRS}
\end{small}
\\
\begin{small}
\textsl{INSA-Lyon, ICJ, UMR5208, F-69621, France}
\end{small}
}
\date{}

\setmarginsrb{25mm}{15mm}{25mm}{15mm}{0mm}{10mm}{0mm}{10mm}

\newtheorem{theo}{Theorem}[section]
\newtheorem{pr}[theo]{Proposition}

\newtheorem{remk}[theo]{Remark}

\newenvironment{rem}{\begin{remk}\normalfont}{\ \rule{0.5em}{0.5em}\end{remk}}
\newcommand{\dem}{\noindent\textsc{Proof\\}}
\newcommand{\refp}[1]{(\ref{#1})}
\newcommand{\fin}{\ \rule{0.5em}{0.5em}}
\numberwithin{equation}{section}

\renewcommand{\a}{\alpha}
\renewcommand{\b}{\beta}
\newcommand{\g}{\gamma}

\renewcommand{\d}{\delta}

\newcommand{\e}{\varepsilon}

\newcommand{\ind}{1\hspace{-.27em}\mbox{\rm l}}

\newcommand{\dis}{\displaystyle}
\newcommand{\lqn}[1]{\noalign{\noindent $\displaystyle{#1}$}}

\begin{document}
\maketitle
\begin{abstract}
We give a method for computing the iterated Laplace transform of the sojourn
time in an union of intervals for linear diffusion processes. This random
variable comes from a model occurring in biology concerning the clustering
of membrane receptors. The way used hinges on solving differential equations.
We finally have a look on the particular case of Brownian motion and we provide
a representation for the Laplace transform of its local time in a finite set.
\end{abstract}

\begin{footnotesize}\sc
\noindent AMS 2000 subject classifications:
\begin{tabular}[t]{l}
primary 60G50; 60J22; \\ secondary 60J10; 60E10.
\end{tabular}  \\
Key words: sojourn time, Laplace transform, linear system.
\end{footnotesize}

\section{Introduction}

We consider a diffusion process $(X_t)_{t\ge 0}$ evolving on the real line $\mathbb{R}$
with infinite lifetime (that is without any killing).
Let $u_1,v_1,\dots,u_n,v_n$ be real numbers such that $u_1<v_1<\dots<u_n<v_n$.
Set also $v_0=-\infty$ and $u_{n+1}=+\infty$. Let us now introduce the union
of intervals
$$
E=\bigcup_{i=1}^n [u_i,v_i]
$$
together with the sojourn time that the process $(X_t)_{t\ge 0}$ spends in~$E$
up to a fixed time~$t$:
$$
T_t=\int_0^t \ind_{E}(X_s)\,ds.
$$
Of course, introducing the sojourn times of $(X_t)_{t\ge 0}$
in $[u_i,v_i]$, $i\in\{1,\dots,n\}$,
$$
T_t^i=\int_0^t \ind_{[u_i,v_i]}(X_s)\,ds,
$$
it is plain that time $T_t$ can be decomposed into
$$
T_t=\sum_{i=1}^nT_t^i.
$$
The functional $T_t$ is one of the most classical and historical random variables
introduced in stochastic processes theory. It has been of great interest for
many researchers in particular when the set~$E$ is an infinite interval
$[v,+\infty)$ or $(-\infty,u]$--cases leading to the famous arcsine law when
$(X_t)_{t\ge 0}$ is Brownian motion for instance--or a bounded interval $[u,v]$,
or when $E$ is reduced to a point $\{u\}$. In this last case, it is related
to the famous local time after some normalization.

The case of an union of intervals for $E$ has appeared to us in a biological
context that we describe here. Consider a cellular medium. The plasmic membrane is a place where
many interactions occur between the cell and its direct external medium.
In classical models, the membrane is described as a fluid mosaic made of
different constituents; certain, the so-called \textsl{ligands} which
possibly come from the extracellular environment, randomly roam
along the membrane while other group into still clusters (the so-called
\textsl{rafts} viewed as \textsl{receptors}) which form inhomogeneities on the
membrane. Few functional properties of these inhomogeneities are known.
A particular phenomenon is the binding mechanism between
ligands and receptors. As a matter of fact, a quantity of great interest for biologists
is the \textsl{docking-time}, this is the time that ligands and
receptors bind up to a finite time (the duration of the experience).
Indeed, this quantity is an important indicator of affinity/sensitivity of
the ligands for receptors. For a more accurate description of this biological context,
we refer the reader to the paper by Car\'e \& Soula~\cite{hedi}.

In~\cite{cholesterol}, we proposed a random walk model for simulating
this problem: the membrane is viewed as the integer line $\mathbb{Z}$
(or a finite torus $\mathbb{Z}/N\mathbb{Z}$), the receptors are fixed numbers
$a_1,\dots,a_r$ in $\mathbb{Z}$ and the ligands evolve along the membrane like
independent random walks $(S_1(i))_{i\in\mathbb{N}},(S_2(i))_{i\in\mathbb{N}},\dots$
on $\mathbb{Z}$. Hence, the movement of the family of $\ell$ ligands is modeled by the
$\ell$-dimensional random walk $(\mathbf{S}(i))_{i\in\mathbb{N}}$ defined,
for any $i\in\mathbb{N}$, by $\mathbf{S}(i)=(S_1(i),\dots,S_\ell(i))$.
The docking-time is proportional to the sojourn time spent by
the random walk $(\mathbf{S}(i))_{i\in\mathbb{N}}$ in the set
$\cup_{j=1}^{\ell} \mathbb{Z}^{j-1}\times\{a_1,\dots,a_r\}
\times\mathbb{Z}^{\ell-j}$. We provided a methodology for computing the
probability distribution of this sojourn time. At the end of~\cite{cholesterol},
we addressed the continuous counterpart to this problem: replace the $\ell$-dimensional
random walk by an $\ell$-dimensional Brownian motion (or a more general
diffusion process) and compute the distribution
of the time spent by Brownian motion in a set of the form
$\cup_{j=1}^{\ell} \mathbb{R}^{j-1}\times E \times\mathbb{R}^{\ell-j}$
where $E$ is a set modeling rafts, that is, $E$ is an union of
intervals.

In this paper, we consider the case of the dimension $\ell=1$;
this is the evolution of one ligand roaming as a diffusion process along
a linear membrane (viewed as $\mathbb{R}$).
In this case, the docking-time is quantified by the sojourn time
of the diffusion process in an union of intervals, namely $T_t$.
Our aim is to describe the probability distribution of~$T_t$
as well as the joint distribution of~$(T_t,X_t)$. For this, we provide
a representation of the iterated Laplace transforms of these distributions
which could be easily implemented on a computer. Actually, our method
allows us to compute the joint distributions of $(T_t^1,\dots,T_t^n)$
and $(T_t^1,\dots,T_t^n,X_t)$.
We finally have a look on the particular case of Brownian motion and we provide
a representation for the Laplace transform of its local time in a finite set.

\section{Settings}

Let us introduce the transition densities of the process $(X_t)_{t\ge 0}$:
for any $t\ge0$ and $x,y\in\mathbb{R}$,
$$
p(t;x,y)=\mathbb{P}_{\!x}\{X_t\in \mathrm{d} y\}/\mathrm{d} y
=\mathbb{P}\{X_t\in \mathrm{d} y|X_0=x\}/\mathrm{d} y
$$
together with its $\lambda$-potential: for any $\lambda>0$ and $x,y\in\mathbb{R}$,
$$
\rho_{_{\hspace{-.05em}\lambda}}\!(x,y)=\int_0^{\infty} \mathrm{e}^{-\lambda t}\,p(t;x,y)\,\mathrm{d} t.
$$

Let $\mathcal{D}$ be the infinitesimal generator of the diffusion $(X_t)_{t\ge 0}$.
This is a second-order differential operator.
Since the lifetime of $(X_t)_{t\ge 0}$ is infinite, the killing rate vanishes
and then $\mathcal{D}$ has the form
$$
\mathcal{D}_x=\frac12\,\sigma(x)^2 \frac{\partial^2}{\partial x^2}+\tau(x) \frac{\partial}{\partial x}.
$$
The notation $\mathcal{D}_x$ means that the operator $\mathcal{D}$ acts on the
variable $x$. It is well-known that the density $p$ solves the Kolmogorov backward equation
$$
\mathcal{D}_x p(t;x,y)=\frac{\partial p}{\partial t}(t;x,y),
\quad p(0^+;x,y)=\d(x-y)
$$
and that the potential solves
$$
\mathcal{D}_x\rho_{_{\hspace{-.05em}\lambda}}\!(x,y)
=\lambda\rho_{_{\hspace{-.05em}\lambda}}\!(x,y)-\d(x-y).
$$
The potential $\rho_{_{\hspace{-.05em}\lambda}}$ is the Green
function of the operator $\mathcal{D}-\lambda\,\mathrm{id}$.

Let us introduce a basis of fundamental solutions $\{c,d\}$ of
the linear second-order differential equation $\mathcal{D}u=\lambda u$ such that
$c$ is increasing and~$d$ is decreasing.
Then the potential $\rho_{_{\hspace{-.05em}\lambda}}$ admits the following representation
(see, e.g., \cite{borodin} p.~19):
$$
\rho_{_{\hspace{-.05em}\lambda}}\!(x,y)=\frac{1}{w(y)}\,c(\min(x,y))\, d(\max(x,y))
$$
where $w(y)=[c'(y)d(y)-c(y)d'(y)]/\kappa(y)$ with $\kappa(y)=2/\sigma(y)^2$.

We see that the potential fulfills the following properties.
The function $x\mapsto\rho_{_{\hspace{-.05em}\lambda}}(x,y)$ is continuous on~$\mathbb{R}$
and, in particular, at $y$. Moreover, it is differentiable on $\mathbb{R}\!\setminus\!\{y\}$
and its derivative admits a jump at $y$. More precisely, we have the
conditions below at $y$:
\begin{equation}\label{conditions-rho}
\begin{cases}
\dis\rho_{_{\hspace{-.05em}\lambda}}(y^+,y)=\rho_{_{\hspace{-.05em}\lambda}}(y^-,y),
\\[1ex]
\dis\frac{\partial \rho_{_{\hspace{-.05em}\lambda}}}{\partial x}(y^+,y)
-\frac{\partial \rho_{_{\hspace{-.05em}\lambda}}}{\partial x}(y^-,y)=-\kappa(y).
\end{cases}
\end{equation}

We also need to introduce the first hitting time of any level~$u$:
$$
\tau_u=\inf\{t\ge 0:X_t=u\}.
$$
The Laplace transform of the distribution of $\tau_u$ is given by
(\cite{borodin} p.~18)
\begin{equation}\label{TL-tau}
\mathbb{E}_x\!\left(\mathrm{e}^{-\lambda\tau_u}\right)=
\begin{cases}
\dis\frac{c(x)}{c(u)}& \mbox{if $x\le u$},
\\[2ex]
\dis\frac{d(x)}{d(u)}& \mbox{if $x\ge u$}.
\end{cases}
\end{equation}

In this paper we compute the iterated Laplace transform of the random vector $\mathbf{T}_t
=(T_t^1,\dots,T_t^n)$:
$$
\varphi_{_{\hspace{-.05em}\lambda,\boldsymbol{\mu}}}\!(x)
=\int_0^{\infty} \mathrm{e}^{-\lambda t}\,\mathbb{E}_x\!\left(
\mathrm{e}^{-<\boldsymbol{\mu},\mathbf{T}_t>}\right) \mathrm{d} t
$$
as well as that of the random vector $(\mathbf{T}_t,X_t)$:
$$
\psi_{_{\hspace{-.05em}\lambda,\boldsymbol{\mu}}}\!(x,y)
=\int_0^{\infty} \mathrm{e}^{-\lambda t}\left[\mathbb{E}_x\!\left(
\mathrm{e}^{-<\boldsymbol{\mu},\mathbf{T}_t>},X_t\in\mathrm{d} y\right)/\mathrm{d} y\right] \mathrm{d} t.
$$
In the foregoing quantities, $\boldsymbol{\mu}$ is a vectorial argument:
$\boldsymbol{\mu}=(\mu_1,\dots,\mu_n)$ and $<\boldsymbol{\mu},\mathbf{T}_t>$
is the inner product of the vectors $\boldsymbol{\mu}$ and $\mathbf{T}_t$:
$$
<\boldsymbol{\mu},\mathbf{T}_t>=\sum_{i=1}^n \mu_iT_t^i.
$$
Notice that we use bold symbols for vectorial variables while unbold symbols
are related to scalar variables.
The iterated Laplace transform of the random variable $T_t$ can be immediately deduced
from $\varphi_{_{\hspace{-.05em}\lambda,\boldsymbol{\mu}}}\!(x)$ by identifying
all the components of $\boldsymbol{\mu}$: if $\mu_1=\dots=\mu_n=\mu$, then
$<\boldsymbol{\mu},\mathbf{T}_t>=\mu T_t.$

For this, we begin by writing integral equations satisfied by the functions
$\varphi_{_{\hspace{-.05em}\lambda,\boldsymbol{\mu}}}$ and
$\psi_{_{\hspace{-.05em}\lambda,\boldsymbol{\mu}}}$ which lead to ordinary linear
differential equations associated with some regularity conditions.
Of course, we have
$$
\varphi_{_{\hspace{-.05em}\lambda,\boldsymbol{\mu}}}\!(x)
=\int_{-\infty}^{+\infty} \psi_{_{\hspace{-.05em}\lambda,\boldsymbol{\mu}}}\!(x,y) \,\mathrm{d} y.
$$
So, because of this relationship, it should be enough to compute
$\psi_{_{\hspace{-.05em}\lambda,\boldsymbol{\mu}}}$
and next integrate it with respect to $y$ for deriving
$\varphi_{_{\hspace{-.05em}\lambda,\boldsymbol{\mu}}}$.
Actually, this way appears to be untractable from the obtained formulas,
so we shall evaluate separately
both functions $\varphi_{_{\hspace{-.05em}\lambda,\boldsymbol{\mu}}}$ and
$\psi_{_{\hspace{-.05em}\lambda,\boldsymbol{\mu}}}$.

\section{Some equations}

In this section, we derive some integral equations (Subsection~\ref{section-int})
and next deduce some differential equations (Subsection~\ref{section-diff}) satisfied
by the functions $\varphi_{_{\hspace{-.05em}\lambda,\boldsymbol{\mu}}}$ and
$\psi_{_{\hspace{-.05em}\lambda,\boldsymbol{\mu}}}$. The method we use is standard, this
is the famous Feynman-Kac approach. Nonetheless, in order to facilitate the
reading of the paper, we provide some details. Our results come from solving
differential equations~\refp{equation-diff}, \refp{equation-diff-bis} below
associated with the boundary values~\refp{conditions}, \refp{conditions-bis},
\refp{conditions-ter}. The aforementioned equations as well as the additional
conditions directly come from integral equations~\refp{eq-phi} and~\refp{eq-psi}.

\subsection{Integral equations}\label{section-int}

%
\begin{pr}
The functions $\varphi_{_{\hspace{-.05em}\lambda,\boldsymbol{\mu}}}$ and
$\psi_{_{\hspace{-.05em}\lambda,\boldsymbol{\mu}}}$ solve the following integral equations:
\begin{align}
\varphi_{_{\hspace{-.05em}\lambda,\boldsymbol{\mu}}}\!(x)
&
=\frac{1}{\lambda}-\sum_{i=1}^n\mu_i \int_{u_i}^{v_i} \rho_{_{\hspace{-.05em}\lambda}}\!(x,z)
\varphi_{_{\hspace{-.05em}\lambda,\boldsymbol{\mu}}}\!(z)\,\mathrm{d} z,
\label{eq-phi}\\
\psi_{_{\hspace{-.05em}\lambda,\boldsymbol{\mu}}}\!(x,y)
&
=\rho_{_{\hspace{-.05em}\lambda}}\!(x,y)-\sum_{i=1}^n\mu_i
\int_{u_i}^{v_i} \rho_{_{\hspace{-.05em}\lambda}}\!(x,z)
\psi_{_{\hspace{-.05em}\lambda,\boldsymbol{\mu}}}\!(z,y)\,\mathrm{d} z.
\label{eq-psi}
\end{align}
\end{pr}
%
\dem
Let us introduce the quantity
$$
\chi_{_{\hspace{-.05em}\boldsymbol{\mu}}}\!(t;x,y)=
\mathbb{E}_x\!\left(\mathrm{e}^{-<\boldsymbol{\mu},\mathbf{T}_t>},
X_t\in\mathrm{d} y\right)\!/\mathrm{d} y.
$$
We write that
\begin{align*}
1-\mathrm{e}^{-<\boldsymbol{\mu},\mathbf{T}_t>}
&
= \int_0^t \left(\sum_{i=1}^n \mu_i\ind_{[u_i,v_i]}(X_s)\right)
\exp\!\left(-\sum_{i=1}^n \mu_i \int_s^t \ind_{[u_i,v_i]}(X_u)\,\mathrm{d}u\right)\mathrm{d}s
\\
&
=\int_0^t \left(\sum_{i=1}^n \mu_i\ind_{[u_i,v_i]}(X_s)\right)
\mathrm{e}^{-<\boldsymbol{\mu},\mathbf{T}_t-\mathbf{T}_s>}\,\mathrm{d}s
=\int_0^t \left(\sum_{i=1}^n \mu_i\ind_{[u_i,v_i]}(X_s)\right)
\,\mathrm{e}^{-<\boldsymbol{\mu},\mathbf{T}_{t-s}\circ \theta_s>}\,\mathrm{d}s,
\end{align*}
where $(\theta_s)_{s\ge 0}$ is the usual shift operator defined by
$X_t\circ \theta_s=X_{s+t}$ for any $s,t\ge 0$. Next, we apply the Markov property:
\begin{align}
\chi_{_{\hspace{-.05em}\boldsymbol{\mu}}}\!(t;x,y)-p(t;x,y)
&
=-\mathbb{E}_x\!\left(1-\mathrm{e}^{-<\boldsymbol{\mu},\mathbf{T}_t>},
X_t\in\mathrm{d}y\right)\!/\mathrm{d}y
\nonumber\\
&
=-\sum_{i=1}^n \mu_i \int_0^t \left[\mathbb{E}_x\!\left(\ind_{[u_i,v_i]}(X_s)\,
\mathrm{e}^{-<\boldsymbol{\mu},\mathbf{T}_{t-s}\circ \theta_s>},
X_{t-s}\circ \theta_s\in\mathrm{d}y\right)\!/\mathrm{d}y\right]\mathrm{d}s
\nonumber\\
&
=-\sum_{i=1}^n \mu_i \int_0^t\!\!\int_{u_i}^{v_i} \mathbb{P}_{\!x}\{X_s\in\mathrm{d}z\}
\left[\mathbb{E}_z\!\left(\mathrm{e}^{-<\boldsymbol{\mu},\mathbf{T}_{t-s}>},
X_{t-s}\in\mathrm{d}y\right)\!/\mathrm{d}y\right]\mathrm{d}s
\nonumber\\
&
=-\sum_{i=1}^n \mu_i \int_0^t\!\!\int_{u_i}^{v_i} p(s;x,z)\,
\chi_{_{\hspace{-.05em}\boldsymbol{\mu}}}\!(t-s;z,y)\,\mathrm{d}s\,\mathrm{d}z.
\label{equation-integral}
\end{align}
By applying the Laplace transform to~\refp{equation-integral},
we get Eq.~\refp{eq-psi} which in turn
yields Eq.~\refp{eq-phi} by integrating with respect to~$y$ on~$\mathbb{R}$ and using
$$
\int_{-\infty}^{+\infty} \rho_{_{\hspace{-.05em}\lambda}}\!(x,y) \,\mathrm{d} y
=\int_0^{\infty} \mathrm{e}^{-\lambda t}\,\mathbb{P}_{\!x}\{X_t\in \mathbb{R}\}\mathrm{d} t=\frac{1}{\lambda}.
$$
\fin

We also state the following result which will be useful.
%
\begin{pr}
For $x\in(-\infty,u_1)$, we have
\begin{align}
\varphi_{_{\hspace{-.05em}\lambda,\boldsymbol{\mu}}}\!(x)
&
=\mathbb{E}_x\!\left(\mathrm{e}^{-\lambda\tau_{u_1}}\right)\!
\Big(\varphi_{_{\hspace{-.05em}\lambda,\boldsymbol{\mu}}}\!(u_1)-\frac{1}{\lambda}\Big)
+\frac{1}{\lambda},
\label{rel-phi}\\
\psi_{_{\hspace{-.05em}\lambda,\boldsymbol{\mu}}}\!(x,y)
-\rho_{_{\hspace{-.05em}\lambda}}\!(x,y)
&
=\mathbb{E}_x\!\left(\mathrm{e}^{-\lambda\tau_{u_1}}\right)\!
\left[\psi_{_{\hspace{-.05em}\lambda,\boldsymbol{\mu}}}\!(u_1,y)
-\rho_{_{\hspace{-.05em}\lambda}}\!(u_1,y)\right]\!.
\label{rel-psi}
\end{align}
\end{pr}
%
\dem
Suppose $x<u_1$. When $\tau_{u_1}>t$, the process $(X_s)_{s\ge 0}$ remains confined
in $(-\infty,u_1)$ up to time $t$; we have $T_t^i=0$ for all
$i\in\{1,\dots,n\}$ and then $1-\mathrm{e}^{-<\boldsymbol{\mu},\mathbf{T}_t>}=0$.
So, we can write that
\begin{align*}
\mathbb{E}_x \!\left(1-\mathrm{e}^{-<\boldsymbol{\mu},\mathbf{T}_t>},
X_t\in\mathrm{d}y\right)\!/\mathrm{d}y
&
=\mathbb{E}_x \!\left(1-\mathrm{e}^{-<\boldsymbol{\mu},\mathbf{T}_t>},
X_t\in\mathrm{d}y,\tau_{u_1}\le t\right)\!/\mathrm{d}y
\\
&
=\mathbb{E}_x \!\left(1-\mathrm{e}^{-<\boldsymbol{\mu},\mathbf{T}_{t-\tau_{u_1}}\circ \theta_{\tau_{u_1}}>},
X_{t-\tau_{u_1}}\circ \theta_{\tau_{u_1}}\in\mathrm{d}y,\tau_{u_1}\le t\right)\!/\mathrm{d}y
\\
&
=\int_0^t \mathbb{P}_{\!x}\{\tau_{u_1}\in\mathrm{d} s\}\,\mathbb{E}_{u_1}
\!\left(1-\mathrm{e}^{-<\boldsymbol{\mu},\mathbf{T}_{t-s}>},
X_{t-s}\in\mathrm{d}y\right)\!/\mathrm{d}y.
\end{align*}
Taking the Laplace transform with respect to $t$, this yields
\begin{align*}
\rho_{_{\hspace{-.05em}\lambda}}\!(x,y)
-\psi_{_{\hspace{-.05em}\lambda,\boldsymbol{\mu}}}\!(x,y)
&
=\int_0^{\infty} \mathrm{e}^{-\lambda t}\left[\mathbb{E}_x
\!\left(1-\mathrm{e}^{-<\boldsymbol{\mu},\mathbf{T}_t>},
X_t\in\mathrm{d}y\right)\!/\mathrm{d}y\right] \mathrm{d} t
\\
&
=\int_0^{\infty} \mathrm{e}^{-\lambda t}\, \mathbb{P}_{\!x}\{\tau_{u_1}\in\mathrm{d} t\}
\int_0^{\infty} \mathrm{e}^{-\lambda t} \left[\mathbb{E}_{u_1}
\!\left(1-\mathrm{e}^{-<\boldsymbol{\mu},\mathbf{T}_t>},
X_t\in\mathrm{d}y\right)\!/\mathrm{d}y\right]\mathrm{d} t
\\
&
=\mathbb{E}_x\!\left(\mathrm{e}^{-\lambda\tau_{u_1}}\right)\!
\left[\rho_{_{\hspace{-.05em}\lambda}}\!(u_1,y)
-\psi_{_{\hspace{-.05em}\lambda,\boldsymbol{\mu}}}\!(u_1,y)\right]\!.
\end{align*}
This proves~\refp{rel-psi} and \refp{rel-phi} ensues by integrating
\refp{rel-psi} with respect to $y$ on $\mathbb{R}$.
\fin

\subsection{Differential equations}\label{section-diff}

We deduce from Eqs~\refp{eq-phi} and~\refp{eq-psi} the famous backward
Kolmogorov equations.
%
\begin{pr}
The functions $\varphi_{_{\hspace{-.05em}\lambda,\boldsymbol{\mu}}}$ and
$\psi_{_{\hspace{-.05em}\lambda,\boldsymbol{\mu}}}$ solve the following differential equations:
\begin{align*}
\mathcal{D}_x\varphi_{_{\hspace{-.05em}\lambda,\boldsymbol{\mu}}}(x)
&
=\bigg[\lambda+\sum_{i=1}^n\mu_i\ind_{[u_i,v_i]}(x)\bigg]
\varphi_{_{\hspace{-.05em}\lambda,\boldsymbol{\mu}}}\!(x)-1,
\\
\mathcal{D}_x\psi_{_{\hspace{-.05em}\lambda,\boldsymbol{\mu}}}(x,y)
&
=\bigg[\lambda+\sum_{i=1}^n\mu_i\ind_{[u_i,v_i]}(x)\bigg]
\psi_{_{\hspace{-.05em}\lambda,\boldsymbol{\mu}}}\!(x,y)-\d_y(x).
\end{align*}
\end{pr}
%
These equations read more explicitly as follows: for $x\in\mathbb{R}$,
\begin{equation}\label{equation-diff}
\mathcal{D}_x\varphi_{_{\hspace{-.05em}\lambda,\boldsymbol{\mu}}}(x)
=\begin{cases}
(\lambda+\mu_i)\,\varphi_{_{\hspace{-.05em}\lambda,\boldsymbol{\mu}}}\!(x)-1
&\mbox{if $x\in(u_i,v_i)$ and $i\in\{1,\dots,n\}$,}
\\
\lambda\,\varphi_{_{\hspace{-.05em}\lambda,\boldsymbol{\mu}}}\!(x)-1
&\mbox{if $x\in(v_i,u_{i+1})$ and $i\in\{0,1,\dots,n\}$,}
\end{cases}
\end{equation}
and for $x\in\mathbb{R}\!\setminus\!\{y\}$,
\begin{equation}\label{equation-diff-bis}
\mathcal{D}_x\psi_{_{\hspace{-.05em}\lambda,\boldsymbol{\mu}}}(x,y)
=\begin{cases}
(\lambda+\mu_i)\,\psi_{_{\hspace{-.05em}\lambda,\boldsymbol{\mu}}}\!(x,y)
&\mbox{if $x\in(u_i,v_i)$ and $i\in\{1,\dots,n\}$,}
\\
\lambda\,\psi_{_{\hspace{-.05em}\lambda,\boldsymbol{\mu}}}\!(x,y)
&\mbox{if $x\in(v_i,u_{i+1})$ and $i\in\{0,1,\dots,n\}$.}
\end{cases}
\end{equation}

Additionally, we see by~\refp{eq-phi} and the regularity properties of
$\rho_{_{\hspace{-.05em}\lambda}}$ that
$\varphi_{_{\hspace{-.05em}\lambda,\boldsymbol{\mu}}}$ is differentiable on~$\mathbb{R}$.
So, we get the following conditions at points $u_i,v_i$, $i\in\{1,\dots,n\}$:
\begin{equation}\label{conditions}
\begin{cases}
\varphi_{_{\hspace{-.05em}\lambda,\boldsymbol{\mu}}}\!(u_i^+)
=\varphi_{_{\hspace{-.05em}\lambda,\boldsymbol{\mu}}}\!(u_i^-),\quad
\varphi_{_{\hspace{-.05em}\lambda,\boldsymbol{\mu}}}'\!(u_i^+)
=\varphi_{_{\hspace{-.05em}\lambda,\boldsymbol{\mu}}}'\!(u_i^-),
\\[.5ex]
\varphi_{_{\hspace{-.05em}\lambda,\boldsymbol{\mu}}}\!(v_i^+)
=\varphi_{_{\hspace{-.05em}\lambda,\boldsymbol{\mu}}}\!(v_i^-),\quad
\varphi_{_{\hspace{-.05em}\lambda,\boldsymbol{\mu}}}'\!(v_i^+)
=\varphi_{_{\hspace{-.05em}\lambda,\boldsymbol{\mu}}}'\!(v_i^-).
\end{cases}
\end{equation}
In the same way, we see by~\refp{eq-psi} that the function
$x\mapsto\psi_{_{\hspace{-.05em}\lambda,\boldsymbol{\mu}}}(x,y)$ is
differentiable on $\mathbb{R}\!\setminus\!\{y\}$.
So, we get the following conditions at points $u_i,v_i$, $i\in\{1,\dots,n\}$:
\begin{equation}\label{conditions-bis}
\begin{cases}
\psi_{_{\hspace{-.05em}\lambda,\boldsymbol{\mu}}}\!(u_i^+,y)
=\psi_{_{\hspace{-.05em}\lambda,\boldsymbol{\mu}}}\!(u_i^-,y),\quad
\psi_{_{\hspace{-.05em}\lambda,\boldsymbol{\mu}}}'\!(u_i^+,y)
=\psi_{_{\hspace{-.05em}\lambda,\boldsymbol{\mu}}}'\!(u_i^-,y),
\\[.5ex]
\psi_{_{\hspace{-.05em}\lambda,\boldsymbol{\mu}}}\!(v_i^+,y)
=\psi_{_{\hspace{-.05em}\lambda,\boldsymbol{\mu}}}\!(v_i^-,y),\quad
\psi_{_{\hspace{-.05em}\lambda,\boldsymbol{\mu}}}'\!(v_i^+,y)
=\psi_{_{\hspace{-.05em}\lambda,\boldsymbol{\mu}}}'\!(v_i^-,y).
\end{cases}
\end{equation}
Moreover, due to~\refp{conditions-rho} and~\refp{eq-psi},
we have the following conditions at point $y$:
\begin{equation}\label{conditions-ter}
\begin{cases}
\psi_{_{\hspace{-.05em}\lambda,\boldsymbol{\mu}}}(y^+,y)=
\psi_{_{\hspace{-.05em}\lambda,\boldsymbol{\mu}}}(y^-,y),
\\[1ex]
\dis\frac{\partial \psi_{_{\hspace{-.05em}\lambda,\boldsymbol{\mu}}}}{\partial x}(y^+,y)
-\frac{\partial \psi_{_{\hspace{-.05em}\lambda,\boldsymbol{\mu}}}}{\partial x}(y^-,y)=-\kappa(y).
\end{cases}
\end{equation}
Different cases for $y$ must be distinguished: $y\in[u_{i_0},v_{i_0}]$
for a certain $i_0\in\{1,\dots,n\}$ or $y\in(v_{i_0},u_{i_0+1})$ for a
certain $i_0\in\{0,1,\dots,n\}$.

In the forthcoming sections, we solve Eq.~\refp{equation-diff} with
conditions~\refp{conditions} and  Eq.~\refp{equation-diff-bis} with
conditions~\refp{conditions-bis}--\refp{conditions-ter}. Let us introduce a basis
of fundamental solutions $\{a_i,b_i\}$ of the linear second-order differential equation
$\mathcal{D}u(x)=(\lambda+\mu_i)\,u(x)$ for any $i\in\{1,\dots,n\}$
as well as a basis of fundamental solutions
$\{c,d\}$ of the equation $\mathcal{D}u(x)=\lambda u(x)$. These functions
are chosen such that $a_i,c$ are increasing and $b_i,d$ are decreasing.

\section{Solving differential equations \refp{equation-diff} and \refp{equation-diff-bis}}

\subsection{Solving Eq.~\refp{equation-diff} with conditions~\refp{conditions}}

In view of the differential operator $\mathcal{D}$, we see that the form of
the solution of~\refp{equation-diff} is
$$
\varphi_{_{\hspace{-.05em}\lambda,\boldsymbol{\mu}}}\!(x)
=\begin{cases}
\dis \a_i \,a_i(x)+\b_i \,b_i(x) +\frac{1}{\lambda+\mu_i} &\mbox{for $x\in (u_i,v_i)$ and $i\in\{1,\dots,n\}$,}
\\[2ex]
\dis \g_i \,c(x)+\d_i \,d(x) +\frac{1}{\lambda} &\mbox{for $x\in (v_i,u_{i+1})$ and $i\in\{0,1,\dots,n\}$.}
\end{cases}
$$
We have to determine the unknown coefficients
$\g_0,\d_0,\a_1,\b_1,\g_1,\d_1,\dots,\a_n,\b_n,\g_n,\d_n$.
Put $\nu_i=\frac{\mu_i}{\lambda(\lambda+\mu_i)}.$

For large enough negative $x$ (so that $x<u_1$), we have
$\mathbb{E}_x\!\left(\mathrm{e}^{-\lambda\tau_{u_1}}\right)=c(x)/c(u_1)$
and \refp{rel-phi} supplies
$$
\varphi_{_{\hspace{-.05em}\lambda,\boldsymbol{\mu}}}\!(x)
=\Big(\frac{1}{\lambda}-\varphi_{_{\hspace{-.05em}\lambda,\boldsymbol{\mu}}}\!(u_1)\Big)
\frac{c(x)}{c(u_1)}+\frac{1}{\lambda}.
$$
This clearly implies that $\d_0=0$. Similarly, considering
$\varphi_{_{\hspace{-.05em}\lambda,\boldsymbol{\mu}}}(x)$ for large positive $x$, we see that
$\g_n=0$.

Next the regularity
conditions~\refp{conditions} at~$u_i$ and~$v_i$ yield
\begin{equation}\label{system}
\begin{cases}
\a_i \,a_i(u_i)+\b_i \,b_i(u_i)-\g_{i-1} \,c(u_i)-\d_{i-1} \,d(u_i)
&\!\!\!\!=\nu_i,
\\
\a_i \,a_i'(u_i)+\b_i \,b_i'(u_i)-\g_{i-1} \,c'(u_i)-\d_{i-1} \,d'(u_i)
&\!\!\!\!=0,
\\
\a_i \,a_i(v_i)+\b_i \,b_i(v_i)-\g_i \,c(v_i)-\d_i \,d(v_i)
&\!\!\!\!=\nu_i,
\\
\a_i \,a_i'(v_i)+\b_i \,b_i'(v_i)-\g_i \,c'(v_i)-\d_i \,d'(v_i)
&\!\!\!\!=0.
\end{cases}
\end{equation}
Let us introduce the matrices
$$
A_i=\begin{pmatrix} \a_i \\ \b_i \end{pmatrix}\!,\quad
B_i=\begin{pmatrix} \g_i \\ \d_i \end{pmatrix}\!,\quad
C_0=\begin{pmatrix} 1 \\ 0 \end{pmatrix}\!,\quad
M_i(x)=\begin{pmatrix} a_i(x) & b_i(x) \\ a_i'(x) & b_i'(x) \end{pmatrix}\!,\quad
N(x)=\begin{pmatrix} c(x) & d(x) \\ c'(x) & d'(x) \end{pmatrix}\!.
$$
Notice that $M_i(x)$ depends on $\lambda+\mu_i$ and $N(x)$ depends on~$\lambda$.
The system~\refp{system} can be rewritten into a matrix form as
\begin{equation}\label{matrix-equation}
\begin{cases}
M_i(u_i)A_i-N(u_i)B_{i-1}=\nu_i C_0,
\\[.5ex]
M_i(v_i)A_i-N(v_i)B_i=\nu_i C_0.
\end{cases}
\end{equation}
We extract from the first equation of~\refp{matrix-equation}
the  relationship $A_i=M_i(u_i)^{-1}N(u_i)B_{i-1}+\nu_i M_i(u_i)^{-1}C_0$.
Plugging this equality into the second equation of~\refp{matrix-equation},
we derive the recursive identity for $B_i$
\begin{equation}\label{recurrence}
B_i=[N(v_i)^{-1}M_i(v_i)M_i(u_i)^{-1}N(u_i)]B_{i-1}+\nu_i[N(v_i)^{-1}M_i(v_i)M_i(u_i)^{-1}-N(v_i)^{-1}]C_0.
\end{equation}
Setting
\begin{equation}\label{def-Pi-Qi}
P_i=N(v_i)^{-1}M_i(v_i)M_i(u_i)^{-1}N(u_i),\quad
Q_i=N(v_i)^{-1}M_i(v_i)M_i(u_i)^{-1}-N(v_i)^{-1},
\end{equation}
we write~\refp{recurrence} more concisely as
\begin{equation}\label{matrix-recurrence}
B_i=P_i B_{i-1}+\nu_i Q_iC_0.
\end{equation}
By iterating this recursive equality, we get for $i\in\{0,1,\dots,n\}$
$$
B_i=P_iP_{i-1}\dots P_1 B_0+(\nu_iQ_i+\nu_{i-1}P_iQ_{i-1}+\nu_{i-2}P_iP_{i-1}Q_{i-2}
+\dots +\nu_1P_iP_{i-1}\dots P_2Q_1)C_0
$$
with the conventions that $P_iP_{i-1}\dots P_2=O$ and $P_iP_{i-1}\dots P_1=I$ if $i=0$,
and $P_iP_{i-1}\dots P_2=I$ if $i=1$. The condition $\d_0=0$ gives $B_0=\g_0 C_0$.
Put, for $i\in\{0,1,\dots,n\}$,
\begin{equation}\label{def-Ri-Si}
R_i=P_iP_{i-1}\dots P_1,\quad
S_i=\nu_iQ_i+\nu_{i-1}P_iQ_{i-1}+\nu_{i-2}P_iP_{i-1}Q_{i-2}+\dots +\nu_1P_iP_{i-1}\dots P_2Q_1.
\end{equation}
By the aforementioned conventions, we have $R_0=I$, $S_0=O$, $S_1=\nu_1Q_1$.
With these settings at hand, we get an simple representation for $B_i$:
$$
B_i=[\g_0 R_i+S_i]C_0.
$$
Now, we use the condition $\g_n=0$. By observing that
$$
\g_n=\big(1\;\;0\big)B_n=\g_0\big(1\;\;0\big)R_n C_0+\big(1\;\;0\big)S_n C_0,
$$
we deduce the value of the coefficient $\g_0$:
$$
\g_0=-\frac{\big(1\;\;0\big)S_nC_0}{\big(1\;\;0\big)R_n C_0}.
$$
Let us point out that $\big(1\;\;0\big)R_nC_0$ (resp. $\big(1\;\;0\big)S_nC_0$) is nothing but
the first entry of the matrix $R_n$ (resp. $S_n$).

Finally, we express $A_i$ by means of~$B_i$: due to~\refp{matrix-equation}, we have
\begin{align*}
A_i
&
=M_i(v_i)^{-1}N(v_i)B_i+\nu_i M_i(v_i)^{-1}C_0
\\
&
=\left[\g_0M_i(v_i)^{-1}N(v_i)R_i+M_i(v_i)^{-1}N(v_i)S_i+\nu_i M_i(v_i)^{-1}\right]C_0.
\end{align*}
We sum up the results obtained in this section in the statement below.
%
\begin{theo}
The iterated Laplace transform of the probability distribution of~$\mathbf{T}_t$
is given by
$$
\int_0^{\infty} \mathrm{e}^{-\lambda t}\,
\mathbb{E}_x\!\left(\mathrm{e}^{-<\boldsymbol{\mu},\mathbf{T}_t>}\right) \mathrm{d} t
=\begin{cases}
\dis\big(1\;\;0\big)M_i(x)A_i+\frac{1}{\lambda+\mu_i}
&\mbox{for $x\in [u_i,v_i]$ and $i\in\{1,\dots,n\}$,}
\\[2ex]
\dis\big(1\;\;0\big) N(x)B_i+\frac{1}{\lambda}
&\mbox{for $x\in (v_i,u_{i+1})$ and $i\in\{0,1,\dots,n\}$,}
\end{cases}
$$
with, for $i\in\{1,\dots,n\}$,
$$
A_i=M_i(v_i)^{-1}N(v_i)\bigg[S_i\!\begin{pmatrix}1\\ 0\end{pmatrix}
-\frac{\big(1\;\;0\big)S_n\!\begin{pmatrix}1\\ 0\end{pmatrix}}{
\big(1\;\;0\big)R_n\!\begin{pmatrix}1\\ 0\end{pmatrix}}
R_i\!\begin{pmatrix}1\\ 0\end{pmatrix}\!\!\bigg]
+\frac{\mu_i}{\lambda(\lambda+\mu_i)}M_i(v_i)^{-1}\!\begin{pmatrix}1\\ 0\end{pmatrix}\!,
$$
for $i\in\{0,1,\dots,n\}$,
$$
B_i=S_i\!\begin{pmatrix}1\\0\end{pmatrix}
-\frac{\big(1\;\;0\big)S_n\!\begin{pmatrix}1\\ 0\end{pmatrix}}{
\big(1\;\;0\big)R_n\!\begin{pmatrix}1\\ 0\end{pmatrix}}
R_i\!\begin{pmatrix}1\\0\end{pmatrix}\!,
$$
where the matrices $R_i,S_i$, $i\in\{0,1,\dots,n\}$, are defined
by~\refp{def-Pi-Qi} and~\refp{def-Ri-Si}.
\end{theo}
%
\begin{rem}\label{remark}
The relationships $B_n=[\g_0 R_n+S_n]C_0$ and $B_n=\d_nD_0$ with
$D_0=\begin{pmatrix}0\\1\end{pmatrix}$ imply
$$
\d_n=\begin{pmatrix}0\;\;1\end{pmatrix}\!B_n
=\g_0\begin{pmatrix}0\;\;1\end{pmatrix}\!R_nC_0
+\begin{pmatrix}0\;\;1\end{pmatrix}\!S_nC_0
$$
which gives the following expression of the coefficient $\d_n$:
$$
\d_n=\begin{pmatrix}0\;\;1\end{pmatrix}\!S_nC_0-\frac{\big(1\;\;0\big)S_nC_0}{
\big(1\;\;0\big)R_n C_0} \begin{pmatrix}0\;\;1\end{pmatrix}\!R_nC_0.
$$
Another expression can be derived by reversing the sense of the
algorithm we used for solving the system~\refp{matrix-equation}. Indeed,
let us rewrite~\refp{matrix-recurrence} as $B_{i-1}=P_i^{-1}B_i -\nu_i P_i^{-1}Q_iC_0$.
Set
$$
\tilde{P}_i=P_i^{-1}=N(u_i)^{-1}M_i(u_i)M_i(v_i)^{-1}N(v_i),\quad
\tilde{Q}_i=-P_i^{-1}Q_i=N(u_i)^{-1}M_i(u_i)M_i(v_i)^{-1}-N(u_i)^{-1}.
$$
The matrices $\tilde{P}_i$ and $\tilde{Q}_i$ can be easily deduced from
$P_i$ and $Q_i$ by interchanging $u_i$ and $v_i$. We have the following recurrence:
$$
B_{i-1}=\tilde{P}_iB_i +\nu_i\tilde{Q}_iC_0
$$
which can be successively iterated from $i=1$ to $n$; this yields
$$
B_0=\tilde{R}_nB_n +\tilde{S}_nC_0
$$
where
$$
\tilde{R}_n=\tilde{P}_1\tilde{P}_2\dots \tilde{P}_n,\quad
\tilde{S}_n=\nu_1\tilde{Q}_1+\nu_2\tilde{P}_1\tilde{Q}_2
+\nu_3\tilde{P}_1\tilde{P}_2\tilde{Q}_3+\dots +\nu_n\tilde{P}_1\tilde{P}_2\dots
\tilde{P}_{n-1}\tilde{Q}_n.
$$
Since $B_0=\g_0C_0$ and $B_n=\d_nD_0$, we get
$$
\g_0C_0=\d_n\tilde{R}_nD_0 +\tilde{S}_nC_0
$$
and we extract the simple expression
$$
\d_n=-\frac{\begin{pmatrix}0\;\;1\end{pmatrix}\!\tilde{S}_nC_0}{
\begin{pmatrix}0\;\;1\end{pmatrix}\!\tilde{R}_n D_0}.
$$
\end{rem}

\subsection{Solving Eq.~\refp{equation-diff-bis} with conditions~\refp{conditions-bis}
and~\refp{conditions-ter}}

We have several cases to consider depending on the location of~$y$ inside or outside
the intervals of the set~$E$.

\vspace{\baselineskip}
\noindent\textsl{\textbf{First case:}} $y\in[u_{i_0},v_{i_0}]$ for
a certain $i_0\in\{1,\dots,n\}$
\vspace{\baselineskip}

The form of the solution of~\refp{equation-diff-bis} is
$$
\psi_{_{\hspace{-.05em}\lambda,\boldsymbol{\mu}}}\!(x,y)
=\begin{cases}
\dis \a_i(y) \,a_i(x)+\b_i(y) \,b_i(x) &\mbox{for $x\in (u_i,v_i)$ and $i\in\{1,\dots,n\}\!\setminus\!\{i_0\}$,}
\\
\dis \g_i(y) \,c(x)+\d_i(y) \,d(x) &\mbox{for $x\in (v_i,u_{i+1})$ and $i\in\{0,\dots,n\}$,}
\\
\dis \a_{i_0,1}(y) \,a_{i_0}(x)+\b_{i_0,1}(y) \,b_{i_0}(x) &\mbox{for $x\in (u_{i_0},y)$,}
\\
\dis \a_{i_0,2}(y) \,a_{i_0}(x)+\b_{i_0,2}(y) \,b_{i_0}(x) &\mbox{for $x\in (y,v_{i_0})$.}
\end{cases}
$$
If $y=u_{i_0}$ (resp. $v_{i_0}$), we consider that the interval $(u_{i_0},y)$
(resp. $(y,v_{i_0})$) is empty.
We have to determine the unknown coefficients
$\g_0(y),\d_0(y),\a_1(y),\b_1(y),\g_1(y),\d_1(y),\dots,\a_{i_0,1}(y),\b_{i_0,1}(y),$
$\a_{i_0,2}(y),\b_{i_0,2}(y),\dots,$ $\a_n(y),\b_n(y),\g_n(y),\d_n(y)$.

For large enough negative $x$ (so that $x<\min(y,u_1)$), we have
$\rho_{_{\hspace{-.05em}\lambda}}\!(x,y)=\frac{1}{w(y)}\,c(x)\,d(y)$
and $\mathbb{E}_x\!\left(\mathrm{e}^{-\lambda\tau_{u_1}}\right)=c(x)/c(u_1)$.
Thus \refp{rel-psi} gives
$$
\psi_{_{\hspace{-.05em}\lambda,\boldsymbol{\mu}}}\!(x,y)
=\left[\frac{1}{c(u_1)} \left[\psi_{_{\hspace{-.05em}\lambda,\boldsymbol{\mu}}}\!(u_1,y)
-\rho_{_{\hspace{-.05em}\lambda}}\!(u_1,y)\right]+\frac{1}{w(y)}\,d(y)\right]c(x).
$$
This clearly implies that $\d_0(y)=0$. Similarly, considering
$\psi_{_{\hspace{-.05em}\lambda,\boldsymbol{\mu}}}(x,y)$ for large positive $x$, we see that
$\g_n(y)=0$.

Next, conditions~\refp{conditions-bis} at $u_i,v_i$,
$i\in\{1,\dots,n\}\!\setminus\!\{i_0\}$, yield
$$
\begin{cases}
\a_i(y) \,a_i(u_i)+\b_i(y) \,b_i(u_i)
&\hspace{-.8em}=\g_{i-1}(y) \,c(u_i)+\d_{i-1}(y) \,d(u_i),
\\
\a_i(y) \,a_i'(u_i)+\b_i(y) \,b_i'(u_i)
&\hspace{-.8em}=\g_{i-1}(y) \,c'(u_i)+\d_{i-1}(y) \,d'(u_i),
\\
\a_i(y) \,a_i(v_i)+\b_i(y) \,b_i(v_i)
&\hspace{-.8em}=\g_i(y) \,c(v_i)+\d_i(y) \,d(v_i),
\\
\a_i(y) \,a_i'(v_i)+\b_i(y) \,b_i'(v_i)
&\hspace{-.8em}=\g_i(y) \,c'(v_i)+\d_i(y) \,d'(v_i).
\end{cases}
$$
Similarly, conditions~\refp{conditions-bis} at $u_{i_0},v_{i_0}$ give
$$
\begin{cases}
\a_{i_01}(y) \,a_{i_0}(u_{i_0})+\b_{i_01}(y) \,b_{i_0}(u_{i_0})
&\hspace{-.8em}=\g_{i_0-1}(y) \,c(u_{i_0})+\d_{i_0-1}(y) \,d(u_{i_0}),
\\
\a_{i_01}(y) \,a_{i_0}'(u_{i_0})+\b_{i_01}(y) \,b_{i_0}'(u_{i_0})
&\hspace{-.8em}=\g_{i_0-1}(y) \,c'(u_{i_0})+\d_{i_0-1}(y) \,d'(u_{i_0}),
\\
\a_{i_02}(y) \,a_{i_0}(v_{i_0})+\b_{i_02}(y) \,b_{i_0}(v_{i_0})
&\hspace{-.8em}=\g_{i_0}(y) \,c(v_{i_0})+\d_{i_0}(y) \,d(v_{i_0}),
\\
\a_{i_02}(y) \,a_{i_0}'(v_{i_0})+\b_{i_02}(y) \,b_{i_0}'(v_{i_0})
&\hspace{-.8em}=\g_{i_0}(y) \,c'(v_{i_0})+\d_{i_0}(y) \,d'(v_{i_0}).
\end{cases}
$$
Additionally, conditions~\refp{conditions-ter} at $y$ yield
$$
\begin{cases}
\a_{i_0,1}(y) \,a_{i_0}(y)+\b_{i_0,1}(y) \,b_{i_0}(y)
-\a_{i_0,2}(y) \,a_{i_0}(y)-\b_{i_0,2}(y) \,b_{i_0}(y)
&\hspace{-.8em}=0,
\\[.5ex]
\a_{i_0,1}(y) \,a_{i_0}'(y)+\b_{i_0,1}(y) \,b_{i_0}'(y)
-\a_{i_0,2}(y) \,a_{i_0}'(y)-\b_{i_0,2}(y) \,b_{i_0}'(y)
&\hspace{-.8em}=\kappa(y).
\end{cases}
$$
Using the matrices introduced in the foregoing subsection and setting also
$$
A_{i_01}(y)=\begin{pmatrix} \a_{i_01}(y) \\ \b_{i_01}(y) \end{pmatrix}\!,\quad
A_{i_02}(y)=\begin{pmatrix} \a_{i_02}(y) \\ \b_{i_02}(y) \end{pmatrix}\!,\quad
D_0=\begin{pmatrix} 0 \\ 1 \end{pmatrix}\!,
$$
we rewrite all these equations into a matrix form: for $i\in\{1,\dots,n\}\!\setminus\!\{i_0\}$,
\begin{equation}\label{matrix-equation-bis}
\begin{cases}
M_i(u_i)A_i(y)&\hspace{-.8em}=N(u_i)B_{i-1}(y),
\\[.5ex]
M_i(v_i)A_i(y)&\hspace{-.8em}=N(v_i)B_i(y),
\end{cases}
\end{equation}
for $i=i_0$,
\begin{equation}\label{matrix-equation-ter}
\begin{cases}
M_{i_0}(u_{i_0})A_{i_01}(y)&\hspace{-.8em}=N(u_{i_0})B_{i_0-1}(y),
\\[.5ex]
M_{i_0}(v_{i_0})A_{i_02}(y)&\hspace{-.8em}=N(v_{i_0})B_{i_0}(y),
\end{cases}
\end{equation}
and finally
\begin{equation}\label{matrix-equation-quater}
M_{i_0}(y)A_{i_02}(y)-M_{i_0}(y)A_{i_01}(y)=-\kappa(y)D_0.
\end{equation}
We extract from~\refp{matrix-equation-bis} the  relationships
$B_i(y)=N(v_i)^{-1}M_i(v_i)A_i(y)$ and $A_i(y)=M_i(u_i)^{-1}N(u_i)B_{i-1}(y)$
which entail for $i\in\{1,\dots,n\}\!\setminus\!\{i_0\}$
\begin{equation}\label{recurrence-bis}
B_i(y)=P_i B_{i-1}(y).
\end{equation}
Therefore, since $\d_0(y)=0$ and then $B_0(y)=\g_0(y)C_0$,
$$
B_i(y)=\begin{cases}
\g_0(y)P_iP_{i-1}\dots P_1C_0 & \mbox{if $i\in\{0,1,\dots,i_0-1\}$,}
\\[.5ex]
P_iP_{i-1}\dots P_{i_0+1}B_{i_0}(y) & \mbox{if $i\in\{i_0+1,\dots,n\}$,}
\end{cases}
$$
with, in view of~\refp{matrix-equation-ter} and~\refp{matrix-equation-quater},
\begin{align*}
B_{i_0}(y)
&
=N(v_{i_0})^{-1}M_{i_0}(v_{i_0})A_{i_02}(y)
\\
&
=N(v_{i_0})^{-1}M_{i_0}(v_{i_0}) \!\left[A_{i_01}(y)-\kappa(y)M_{i_0}(y)^{-1}D_0\right]
\\
&
=N(v_{i_0})^{-1}M_{i_0}(v_{i_0}) \!\left[M_{i_0}(u_{i_0})^{-1}N(u_{i_0})B_{i_0-1}(y)
-\kappa(y)M_{i_0}(y)^{-1}D_0\right]
\\
&
=N(v_{i_0})^{-1}M_{i_0}(v_{i_0}) \!\left[\g_0(y)M_{i_0}(u_{i_0})^{-1}N(u_{i_0})
P_{i_0-1}P_{i_0-2}\dots P_1C_0-\kappa(y)M_{i_0}(y)^{-1}D_0\right]
\\
&
=\g_0(y)P_{i_0}P_{i_0-1}\dots P_1C_0
-\kappa(y)N(v_{i_0})^{-1}M_{i_0}(v_{i_0}) M_{i_0}(y)^{-1}D_0.
\end{align*}
As a byproduct, recalling that we set $R_i=P_iP_{i-1}\dots P_1$ and putting,
for $i\in\{0,1,\dots,i_0-1\}$,
\begin{equation}\label{def-Ui}
U_i=P_iP_{i-1}\dots P_{i_0+1}N(v_{i_0})^{-1}M_{i_0}(v_{i_0})
\end{equation}
with the conventions that $P_iP_{i-1}\dots P_{i_0+1}=I$ if $i=i_0$ and
$P_iP_{i-1}\dots P_{i_0+1}=O$ if $i<i_0$, we obtain the following
expression for $B_i(y)$: for any $i\in\{0,1,\dots,n\}$,
\begin{equation}\label{exp-Bi}
B_i(y)=\g_0(y)R_iC_0-\kappa(y)U_i M_{i_0}(y)^{-1}D_0.
\end{equation}
Now it remains to determine $\g_0(y)$. For this, we invoke the condition $\g_n(y)=0$.
Since
$$
\g_n(y)=\big(1\;\;0\big)B_n(y)
=\g_0(y)\big(1\;\;0\big)R_n C_0
-\kappa(y)\big(1\;\;0\big)U_n M_{i_0}(y)^{-1} D_0,
$$
we deduce the value of the coefficient $\g_0(y)$:
$$
\g_0(y)=\kappa(y)\,\frac{\big(1\;\;0\big)U_n M_{i_0}(y)^{-1}D_0}{
\big(1\;\;0\big)R_n C_0}.
$$
Finally, due to~\refp{matrix-equation-bis} and~\refp{exp-Bi},
we can express $A_i(y)$ by means of $B_i(y)$, for
$i\in\{1,\dots,n\}\!\setminus\!\{i_0\}$, as
\begin{align*}
A_i(y)
&
=M_i(v_i)^{-1}N(v_i)B_i(y)
=M_i(v_i)^{-1}N(v_i) \left[\g_0(y)R_iC_0-\kappa(y)U_i M_{i_0}(y)^{-1}D_0\right]\!.
\end{align*}
Moreover, by~\refp{matrix-equation-ter}, \refp{matrix-equation-quater} and~\refp{exp-Bi},
\begin{align*}
A_{i_01}(y)
&
=M_{i_0}(u_{i_0})^{-1}N(u_{i_0})B_{i_0-1}(y)
=\g_0(y)M_{i_0}(u_{i_0})^{-1}N(u_{i_0})R_{i_0-1}C_0,
\\
A_{i_02}(y)
&
=M_{i_0}(v_{i_0})^{-1}N(v_{i_0})B_{i_0}(y)
=\g_0(y)M_{i_0}(u_{i_0})^{-1}N(u_{i_0})R_{i_0-1}C_0-\kappa(y)M_{i_0}(y)^{-1}D_0.
\end{align*}
%
\begin{rem}\label{remark-bis}
As in Remark~\ref{remark}, we can obtain two expressions of the coefficient $\d_n(y)$.
Because of the relationships $B_n(y)=\g_0(y) R_nC_0-\kappa(y)U_nM_{i_0}(y)^{-1}D_0$ and $B_n(y)=\d_n(y)D_0$,
we can see that
$$
\d_n(y)=\kappa(y)\,\frac{\big(1\;\;0\big)U_nM_{i_0}(y)^{-1}D_0}{
\big(1\;\;0\big)R_n C_0} \begin{pmatrix}0\;\;1\end{pmatrix}\!R_nC_0
-\kappa(y)\begin{pmatrix}0\;\;1\end{pmatrix}\!U_nM_{i_0}(y)^{-1}D_0.
$$
Another expression can be derived by reversing the sense of the
algorithm we used for solving the system~\refp{matrix-equation-bis}. Indeed,
let us rewrite~\refp{matrix-recurrence} as $B_{i-1}(y)=\tilde{P}_iB_i(y)$
for $i\in\{1,\dots,n\}\!\setminus\!\{i_0\}$ (recall that $\tilde{P}_i=P_i^{-1}$).
We have
$$
B_0(y)=\tilde{P}_1\dots\tilde{P}_{i_0-1}B_{i_0-1}(y),\quad
B_{i_0}(y)=\tilde{P}_{i_0+1}\dots\tilde{P}_nB_n(y).
$$
Moreover, by~\refp{matrix-equation-ter} and~\refp{matrix-equation-quater},
we successively have
\begin{align*}
B_{i_0-1}(y)
&
=N(u_{i_0})^{-1}M_{i_0}(u_{i_0})A_{i_01}(y)
=N(u_{i_0})^{-1}M_{i_0}(u_{i_0})[A_{i_02}(y)+\kappa(y)M_{i_0}(y)^{-1}D_0]
\\
&
=N(u_{i_0})^{-1}M_{i_0}(u_{i_0})[M_{i_0}(v_{i_0})^{-1}N(v_{i_0})B_{i_0}(y)+\kappa(y)M_{i_0}(y)^{-1}D_0]
\\
&
=\tilde{P}_{i_0}B_{i_0}(y)+\kappa(y)N(u_{i_0})^{-1}M_{i_0}(u_{i_0})M_{i_0}(y)^{-1}D_0
\\
&
=\tilde{P}_{i_0}\tilde{P}_{i_0+1}\dots\tilde{P}_nB_n(y)
+\kappa(y)N(u_{i_0})^{-1}M_{i_0}(u_{i_0})M_{i_0}(y)^{-1}D_0
\end{align*}
and then
$$
B_0(y)=\tilde{P}_1\dots\tilde{P}_nB_n(y)
+\kappa(y)\tilde{P}_1\dots\tilde{P}_{i_0-1}N(u_{i_0})^{-1}M_{i_0}(u_{i_0})M_{i_0}(y)^{-1}D_0.
$$
We rewrite this last equality as
$$
\g_0(y)C_0=\d_n(y)\tilde{R}_n D_0+\kappa(y)\tilde{U}_0M_{i_0}(y)^{-1}D_0
$$
with $\tilde{R}_n=\tilde{P}_1\dots\tilde{P}_n$
and $\tilde{U}_0=\tilde{P}_1\dots\tilde{P}_{i_0-1}N(u_{i_0})^{-1}M_{i_0}(u_{i_0})$.
We finally get the following expression of $\d_n(y)$:
$$
\d_n(y)=-\kappa(y)\,\frac{\begin{pmatrix}0\;\;1\end{pmatrix}\!\tilde{U}_0M_{i_0}(y)^{-1}D_0}{
\begin{pmatrix}0\;\;1\end{pmatrix}\!\tilde{R}_n D_0}.
$$
\end{rem}
%

\vspace{\baselineskip}
\noindent\textsl{\textbf{Second case:}} $y\in(v_{i_0},u_{i_0+1})$ for
a certain $i_0\in\{0,\dots,n\}$
\vspace{\baselineskip}

Since the computations are very analogous, we briefly outline the corresponding matrix equations.
Set
$$
B_{i_01}(y)=\begin{pmatrix} \g_{i_01}(y) \\ \d_{i_01}(y) \end{pmatrix}\!,\quad
B_{i_02}(y)=\begin{pmatrix} \g_{i_02}(y) \\ \d_{i_02}(y) \end{pmatrix}\!.
$$
The equations write as follows:
$$
\begin{cases}
M_i(u_i)A_i(y)=N(u_i)B_{i-1}(y) & \mbox{for $i\in\{1,\dots,n\}\!\setminus\!\{i_0+1\}$,}
\\
M_i(v_i)A_i(y)=N(v_i)B_i(y) & \mbox{for $i\in\{1,\dots,n\}\!\setminus\!\{i_0\}$,}
\\
M_{i_0}(v_{i_0})A_{i_0}(y)=N(v_{i_0})B_{i_01}(y),
\\
M_{i_0+1}(u_{i_0+1})A_{i_0+1}(y)=N(u_{i_0+1})B_{i_02}(y),
\\
N(y)B_{i_02}(y)-N(y)B_{i_01}(y)=-\kappa(y)D_0.
\end{cases}
$$
We deduce that
$$
B_i(y)=\begin{cases}
\g_0(y)P_iP_{i-1}\dots P_1C_0 & \mbox{if $i\in\{0,1,\dots,i_0-1\}$,}
\\[.5ex]
P_iP_{i-1}\dots P_{i_0+2}B_{i_0+1}(y) & \mbox{if $i\in\{i_0+1,\dots,n\}$,}
\end{cases}
$$
with
\begin{align*}
B_{i_0+1}(y)
&
=N(v_{i_0+1})^{-1}M_{i_0+1}(v_{i_0+1})A_{i_0+1}(y)
\\
&
=N(v_{i_0+1})^{-1}M_{i_0+1}(v_{i_0+1})M_{i_0+1}(u_{i_0+1})^{-1}N(u_{i_0+1})B_{i_02}(y)
\\
&
=P_{i_0+1}B_{i_02}(y).
\end{align*}
We need to compute $B_{i_02}(y)$: we have
$$
B_{i_02}(y)=B_{i_01}(y)-\kappa(y)N(y)^{-1}D_0
$$
with
\begin{align*}
B_{i_01}(y)
&
=N(v_{i_0})^{-1}M_{i_0}(v_{i_0})A_{i_0}(y)
=N(v_{i_0})^{-1}M_{i_0}(v_{i_0})M_{i_0}(u_{i_0})^{-1}N(u_{i_0})B_{i_0-1}(y)
\\
&
=P_{i_0}B_{i_0-1}(y)=\g_0(y)P_{i_0}P_{i_0-1}\dots P_1C_0.
\end{align*}
Thus
$$
B_{i_0+1}(y)=P_{i_0+1}\left[\g_0(y)P_{i_0}\dots P_1C_0-\kappa(y)N(y)^{-1}D_0\right]
$$
and for $i\in\{i_0+1,\dots,n\}$,
$$
B_i(y)=\g_0(y)P_iP_{i-1}\dots P_1C_0-\kappa(y)P_iP_{i-1}\dots P_{i_0+1}N(y)^{-1}D_0.
$$
Setting, for $i\in\{0,1,\dots,n\}$,
\begin{equation}\label{def-Vi}
V_i=P_iP_{i-1}\dots P_{i_0+1},
\end{equation}
with the conventions that $V_i=I$ if $i=i_0$ and $V_i=O$ if $i<i_0$,
we obtain the following expression for $B_i(y)$: for any $i\in\{0,1,\dots,n\}$,
$$
B_i(y)=\g_0(y)R_iC_0-\kappa(y)V_i N(y)^{-1}D_0.
$$
The condition $\g_n(y)=0$ supplies the value of the coefficient $\g_0(y)$ within
$B_i(y)$:
$$
\g_0(y)=\kappa(y)\,\frac{\big(1\;\;0\big)V_n N(y)^{-1}D_0}{
\big(1\;\;0\big)R_n C_0}.
$$
%
\begin{rem}
As in Remark~\ref{remark-bis}, we can obtain two expressions of the coefficient $\d_n(y)$.
We only outline the main details. On one hand, the first expression writes
$$
\d_n(y)=\kappa(y)\,\frac{\big(1\;\;0\big)V_nN(y)^{-1}D_0}{
\big(1\;\;0\big)R_n C_0} \begin{pmatrix}0\;\;1\end{pmatrix}\!R_nC_0
-\kappa(y)\begin{pmatrix}0\;\;1\end{pmatrix}\!V_nN(y)^{-1}D_0.
$$
On the other hand, we have
\begin{align*}
B_0(y)
&
=\tilde{P}_1\dots\tilde{P}_{i_0-1}B_{i_0-1}(y)
=\tilde{P}_1\dots\tilde{P}_{i_0}B_{i_01}(y)
=\tilde{P}_1\dots\tilde{P}_{i_0}[B_{i_02}(y)+\kappa(y)N(y)^{-1}D_0]
\\
&
=\tilde{P}_1\dots\tilde{P}_{i_0+1}B_{i_0+1}(y)
+\kappa(y)\tilde{P}_1\dots\tilde{P}_{i_0}N(y)^{-1}D_0
=\tilde{P}_1\dots\tilde{P}_nB_n(y)
+\kappa(y)\tilde{P}_1\dots\tilde{P}_{i_0}N(y)^{-1}D_0
\end{align*}
which we rewrite as
$$
\g_0(y)C_0=\d_n(y)\tilde{R}_n D_0+\kappa(y)\tilde{V}_0N(y)^{-1}D_0
$$
with $\tilde{V}_0=\tilde{P}_1\dots\tilde{P}_{i_0}$.
We get the second expression of $\d_n(y)$:
$$
\d_n(y)=-\kappa(y)\,\frac{\begin{pmatrix}0\;\;1\end{pmatrix}\!\tilde{V}_0N(y)^{-1}D_0}{
\begin{pmatrix}0\;\;1\end{pmatrix}\!\tilde{R}_n D_0}.
$$
\end{rem}
%
Finally, we can express $A_i(y)$ by means of $B_i(y)$: for
$i\in\{1,\dots,n\}\!\setminus\!\{i_0\}$,
$$
A_i(y)=M_i(v_i)^{-1}N(v_i)B_i(y)
=M_i(v_i)^{-1}N(v_i) \left[\g_0(y)R_iC_0-\kappa(y)V_i N(y)^{-1}D_0\right]
$$
and, for $i=i_0$,
$$
A_{i_0}(y)=M_{i_0}(u_{i_0})^{-1}N(u_{i_0})B_{i_0-1}(y)
=M_{i_0}(u_{i_0})^{-1}N(u_{i_0}) \left[\g_0(y)R_{i_0-1}C_0-\kappa(y)V_{i_0-1} N(y)^{-1}D_0\right]\!.
$$
We sum up the results obtained in this section in the following statement.
%
\begin{theo}
The iterated Laplace transform of the joint probability distribution of~$(\mathbf{T}_t,X_t)$
is given by the formulas below.

1) If $y\in[u_{i_0},v_{i_0}]$ for a certain $i_0\in\{1,\dots,n\}$,
$$
\int_0^{\infty} \mathrm{e}^{-\lambda t}\left[\mathbb{E}_x\!
\left(\mathrm{e}^{-<\boldsymbol{\mu},\mathbf{T}_t>},X_t\in\mathrm{d} y\right)
/\mathrm{d} y\right] \mathrm{d} t
=\begin{cases}
\dis\big(1\;\;0\big)M_i(x)A_i(y)
&\mbox{for $x\in [u_i,v_i]$ and $i\in\{1,\dots,n\}\!\setminus\!\{i_0\}$,}
\\
\dis\big(1\;\;0\big) N(x)B_i(y)
&\mbox{for $x\in (v_i,u_{i+1})$ and $i\in\{0,\dots,n\}$,}
\\
\dis\big(1\;\;0\big)M_{i_0}(x)A_{i_01}(y) &\mbox{for $x\in [u_{i_0},y]$,}
\\
\dis\big(1\;\;0\big)M_{i_0}(x)A_{i_02}(y) &\mbox{for $x\in [y,v_{i_0}]$,}
\end{cases}
$$
with, for $i\in\{1,\dots,n\}\!\setminus\!\{i_0\}$,
\begin{align*}
A_i(y)
&
=\kappa(y)M_i(v_i)^{-1}N(v_i) \bigg[
\frac{\big(1\;\;0\big)U_n M_{i_0}(y)^{-1}\!\begin{pmatrix}0\\1\end{pmatrix}}{
\big(1\;\;0\big)R_n\!\begin{pmatrix}1\\0\end{pmatrix}}
R_i\!\begin{pmatrix}1\\0\end{pmatrix}-U_i M_{i_0}(y)^{-1}\!\begin{pmatrix}0\\1\end{pmatrix}\!\!\bigg]\!,
\\
A_{i_01}(y)
&
=\kappa(y)\,\frac{\big(1\;\;0\big)U_n M_{i_0}(y)^{-1}\!\begin{pmatrix}0\\1\end{pmatrix}}{
\big(1\;\;0\big)R_n\!\begin{pmatrix}1\\0\end{pmatrix}}
M_{i_0}(u_{i_0})^{-1}N(u_{i_0})R_{i_0-1}\!\begin{pmatrix}1\\0\end{pmatrix}\!,
\\
A_{i_02}(y)
&
=\kappa(y)\bigg[
\frac{\big(1\;\;0\big)U_n M_{i_0}(y)^{-1}\!\begin{pmatrix}0\\1\end{pmatrix}}{
\big(1\;\;0\big)R_n\!\begin{pmatrix}1\\0\end{pmatrix}}
M_{i_0}(u_{i_0})^{-1}N(u_{i_0})R_{i_0-1}\!\begin{pmatrix}1\\0\end{pmatrix}
-M_{i_0}(y)^{-1}\!\begin{pmatrix}0\\1\end{pmatrix}\!\!\bigg]\!,
\end{align*}
and, for $i\in\{0,1,\dots,n\}$,
$$
B_i(y)=\kappa(y)\bigg[
\frac{\big(1\;\;0\big)U_n M_{i_0}(y)^{-1}\!\begin{pmatrix}0\\1\end{pmatrix}}{
\big(1\;\;0\big)R_n\!\begin{pmatrix}1\\0\end{pmatrix}}
R_i\!\begin{pmatrix}1\\0\end{pmatrix}-U_i M_{i_0}(y)^{-1}\!\begin{pmatrix}0\\1\end{pmatrix}
\!\!\bigg]\!,
$$
where the matrices $R_i,U_i$ are defined by~\refp{def-Pi-Qi}, \refp{def-Ri-Si}
and~\refp{def-Ui}.

2) If $y\in(v_{i_0},u_{i_0+1})$ for a certain $i_0\in\{0,\dots,n\}$,
$$
\int_0^{\infty} \mathrm{e}^{-\lambda t}\left[\mathbb{E}_x\!
\left(\mathrm{e}^{-<\boldsymbol{\mu},\mathbf{T}_t>},X_t\in\mathrm{d} y\right)
/\mathrm{d} y\right] \mathrm{d} t
=\begin{cases}
\dis\big(1\;\;0\big)M_i(x)A_i(y)
&\mbox{for $x\in [u_i,v_i]$ and $i\in\{1,\dots,n\}$,}
\\
\dis\big(1\;\;0\big) N(x)B_i(y)
&\mbox{for $x\in (v_i,u_{i+1})$ and $i\in\{0,\dots,n\}\!\setminus\!\{i_0\}$,}
\\
\dis\big(1\;\;0\big)N(x)B_{i_01}(y) &\mbox{for $x\in (v_{i_0},y]$,}
\\
\dis\big(1\;\;0\big)N(x)B_{i_02}(y) &\mbox{for $x\in [y,u_{i_0+1})$,}
\end{cases}
$$
with, for $i\in\{1,\dots,n\}\!\setminus\!\{i_0\}$,
\begin{align*}
A_i(y)
&
=\kappa(y)M_i(v_i)^{-1}N(v_i) \bigg[
\frac{\big(1\;\;0\big)V_n N(y)^{-1}\!\begin{pmatrix}0\\1\end{pmatrix}}{
\big(1\;\;0\big)R_n\!\begin{pmatrix}1\\0\end{pmatrix}}
R_i\!\begin{pmatrix}1\\0\end{pmatrix}-V_i N(y)^{-1}\!\begin{pmatrix}0\\1\end{pmatrix}\!\!\bigg]\!,
\\
B_i(y)
&
=\kappa(y)\bigg[
\frac{\big(1\;\;0\big)V_n N(y)^{-1}\!\begin{pmatrix}0\\1\end{pmatrix}}{
\big(1\;\;0\big)R_n\!\begin{pmatrix}1\\0\end{pmatrix}}
R_i\!\begin{pmatrix}1\\0\end{pmatrix}-V_i N(y)^{-1}\!\begin{pmatrix}0\\1\end{pmatrix}\!\!\bigg]\!,
\end{align*}
and
\begin{align*}
A_{i_0}(y)
&
=\kappa(y)M_{i_0}(u_{i_0})^{-1}N(u_{i_0}) \bigg[
\frac{\big(1\;\;0\big)V_n N(y)^{-1}\!\begin{pmatrix}0\\1\end{pmatrix}}{
\big(1\;\;0\big)R_n\!\begin{pmatrix}1\\0\end{pmatrix}}
R_{i_0-1}\!\begin{pmatrix}1\\0\end{pmatrix}-V_{i_0-1} N(y)^{-1}\!\begin{pmatrix}0\\1\end{pmatrix}\!\!\bigg],
\\
B_{i_01}(y)
&
=\kappa(y)\,\frac{\big(1\;\;0\big)V_n N(y)^{-1}\!\begin{pmatrix}0\\1\end{pmatrix}}{
\big(1\;\;0\big)R_n\!\begin{pmatrix}1\\0\end{pmatrix}}
R_{i_0}\!\begin{pmatrix}1\\0\end{pmatrix}\!,
\\
B_{i_02}(y)
&
=\kappa(y)\bigg[\frac{\big(1\;\;0\big)V_n N(y)^{-1}\!\begin{pmatrix}0\\1\end{pmatrix}}{
\big(1\;\;0\big)R_n\!\begin{pmatrix}1\\0\end{pmatrix}}
R_{i_0}\!\begin{pmatrix}1\\0\end{pmatrix}-N(y)^{-1}\!\begin{pmatrix}0\\1\end{pmatrix}\!\!\bigg]\!,
\end{align*}
where the matrices $R_i,V_i$ are defined by~\refp{def-Pi-Qi}, \refp{def-Ri-Si}
and~\refp{def-Vi}.
\end{theo}
%

\section{Examples}

\subsection{Case of one bounded interval}

In this part, we focus on the set~$E$ made of one two-sided interval $E=[u,v]$.
This classical case corresponds to the values of the parameters $n=1$ and
$u_1=u$, $v_1=v$, $T_t^1=T_t$. We relabel the argument $\mu_1$ into $\mu$,
the functions $a_1,b_1$ into $a,b$ and the related matrix $M_1$ into $M$.
The settings write here $R_0=I$, $S_0=O$,
$R_1=P_1=N(v)^{-1}M(v)M(u)^{-1}N(u)$, $Q_1=N(v)^{-1}M(v)M(u)^{-1}-N(v)^{-1}$,
$S_1=\frac{\mu}{\lambda(\lambda+\mu)}\,Q_1$,
$\tilde{R}_1=\tilde{P}_1=N(u)^{-1}M(u)M(v)^{-1}N(v)$,
$\tilde{Q}_1=N(u)^{-1}M(u)M(v)^{-1}-N(u)^{-1}$,
$\tilde{S}_1=\frac{\mu}{\lambda(\lambda+\mu)}\,\tilde{Q}_1$.

\vspace{\baselineskip}
\noindent\textsl{\textbf{Probability distribution of $T_t$}}
\vspace{\baselineskip}

The distribution of $T_t$ is characterized by
$$
\int_0^{\infty} \mathrm{e}^{-\lambda t}\, \mathbb{E}_x\!\left(\mathrm{e}^{-\mu T_t}\right) \mathrm{d} t
=\begin{cases}
\dis\big(1\;\;0\big) N(x)B_0+\frac{1}{\lambda} &\mbox{for $x\in (-\infty,u]$,}
\\[2ex]
\dis\big(1\;\;0\big)M(x)A_1+\frac{1}{\lambda+\mu} &\mbox{for $x\in [u,v]$,}
\\[2ex]
\dis\big(1\;\;0\big) N(x)B_1+\frac{1}{\lambda} &\mbox{for $x\in [v,+\infty)$.}
\end{cases}
$$
Observing that $M(v)^{-1}N(v)P_1=M(u)^{-1}N(u)$ and $M(v)^{-1}N(v)Q_1+M(v)^{-1}=M(u)^{-1}$,
the matrix $A_1$ can be simplified into
$$
A_1=\g_0 M(u)^{-1}N(u)\!\begin{pmatrix}1\\0\end{pmatrix}
+\frac{\mu}{\lambda(\lambda+\mu)}M(u)^{-1}\!\begin{pmatrix}1\\0\end{pmatrix}
$$
and the other matrices are given by
$B_0=\g_0\!\begin{pmatrix}1\\0\end{pmatrix}$ and $B_1=\d_1\!\begin{pmatrix}0\\1\end{pmatrix}$
with
$$
\g_0=-\frac{\mu}{\lambda(\lambda+\mu)}\,\frac{\big(1\;\;0\big)Q_1\!\begin{pmatrix}1\\0\end{pmatrix}}{
\big(1\;\;0\big)P_1\!\begin{pmatrix}1\\0\end{pmatrix}},\quad
\d_1=-\frac{\mu}{\lambda(\lambda+\mu)}\,\frac{\begin{pmatrix}0\;\;1\end{pmatrix}\!\tilde{Q}_1\!\begin{pmatrix}1\\0\end{pmatrix}}{
\begin{pmatrix}0\;\;1\end{pmatrix}\!\tilde{P}_1\!\begin{pmatrix}0\\1\end{pmatrix}}.
$$
Moreover,
$$
\big(1\;\;0\big)N(x)B_0=\g_0c(x),\quad\dis\big(1\;\;0\big) N(x)B_1=\d_1d(x).
$$

\vspace{\baselineskip}
\noindent\textsl{\textbf{Joint probability distribution of $(T_t,X_t)$}}
\vspace{\baselineskip}

1) Assume that $y\in(-\infty,u]$. This case corresponds to $i_0=0$ and
we have $V_0=I$, $V_1=P_1$.
The distribution of $(T_t,X_t)$ is characterized by
$$
\int_0^{\infty} \mathrm{e}^{-\lambda t}\left[\mathbb{E}_x\!\left(\mathrm{e}^{-\mu T_t},X_t\in\mathrm{d} y\right)
/\mathrm{d} y\right] \mathrm{d} t
=\begin{cases}
\dis\big(1\;\;0\big) N(x)B_{01}(y) &\mbox{for $x\in (-\infty,y]$,}
\\
\dis\big(1\;\;0\big)N(x)B_{02}(y) &\mbox{for $x\in [y,u]$,}
\\
\dis\big(1\;\;0\big)M(x)A_1(y) &\mbox{for $x\in [u,v]$,}
\\
\dis\big(1\;\;0\big) N(x)B_1(y) &\mbox{for $x\in [v,+\infty)$.}
\end{cases}
$$
Since $M(v)^{-1}N(v)P_1=M(u)^{-1}N(u)$, we have for $A_1(y)$
$$
A_1(y)=\g_0(y)M(v)^{-1}N(v)\!\begin{pmatrix}1\\0\end{pmatrix}
-\kappa(y)M(u)^{-1}N(u) N(y)^{-1}\!\begin{pmatrix}0\\1\end{pmatrix}\!,
$$
and the other matrices are given by
$$
B_{01}(y)=\g_0(y)\!\begin{pmatrix}1\\0\end{pmatrix}\!,\quad
B_{02}(y)=\g_0(y)\!\begin{pmatrix}1\\0\end{pmatrix}-
\kappa(y)N(y)^{-1}\!\begin{pmatrix}0\\1\end{pmatrix}\!,
$$
$$
B_1(y)=P_1\!\left[\g_0(y)\!\begin{pmatrix}1\\0\end{pmatrix}-
\kappa(y)N(y)^{-1}\!\begin{pmatrix}0\\1\end{pmatrix} \!\right]
$$
with
$$
\g_0(y)=\kappa(y)\,\frac{\big(1\;\;0\big)P_1N(y)^{-1}\!\begin{pmatrix}0\\1\end{pmatrix}}{
\big(1\;\;0\big)P_1\!\begin{pmatrix}1\\0\end{pmatrix}}.
$$

2) Assume that $y\in[u,v]$. This case corresponds to $i_0=1$ and we have
$U_0=O$, $U_1=N(v)^{-1}M(v)$.
The distribution of $(T_t,X_t)$ is characterized by
$$
\int_0^{\infty} \mathrm{e}^{-\lambda t}\left[\mathbb{E}_x\!\left(\mathrm{e}^{-\mu T_t},X_t\in\mathrm{d} y\right)
/\mathrm{d} y\right] \mathrm{d} t
=\begin{cases}
\dis\big(1\;\;0\big) N(x)B_0(y) &\mbox{for $x\in (-\infty,u]$,}
\\
\dis\big(1\;\;0\big)M(x)A_{11}(y) &\mbox{for $x\in [u,y]$,}
\\
\dis\big(1\;\;0\big)M(x)A_{12}(y) &\mbox{for $x\in [y,v]$,}
\\
\dis\big(1\;\;0\big) N(x)B_1(y) &\mbox{for $x\in [v,+\infty)$.}
\end{cases}
$$
The matrices are given by
$$
A_{11}(y)=\g_0(y)M(u)^{-1}N(u) \!\begin{pmatrix}1\\0\end{pmatrix}\!,\quad
A_{12}(y)=\g_0(y)M(u)^{-1}N(u) \!\begin{pmatrix}1\\0\end{pmatrix}
-\kappa(y)M(y)^{-1}\!\begin{pmatrix}0\\1\end{pmatrix}\!,
$$
$$
B_0(y)=\g_0(y)\!\begin{pmatrix}1\\0\end{pmatrix}\!,\quad
B_1(y)=\g_0(y)P_1\!\begin{pmatrix}1\\0\end{pmatrix}
-\kappa(y)U_1M(y)^{-1}\!\begin{pmatrix}0\\1\end{pmatrix}
$$
with
$$
\g_0(y)=\kappa(y)\,\frac{\big(1\;\;0\big)U_1M(y)^{-1}\!\begin{pmatrix}0\\1\end{pmatrix}}{
\big(1\;\;0\big)P_1 \!\begin{pmatrix}1\\0\end{pmatrix}}.
$$

3) Assume that $y\in[v,+\infty)$. This case corresponds to $i_0=1$ and we have
$V_0=O$, $V_1=I$. The distribution of $(T_t,X_t)$ is characterized by
$$
\int_0^{\infty} \mathrm{e}^{-\lambda t}\left[\mathbb{E}_x\!\left(\mathrm{e}^{-\mu T_t},X_t\in\mathrm{d} y\right)
/\mathrm{d} y\right] \mathrm{d} t
=\begin{cases}
\dis\big(1\;\;0\big) N(x)B_0(y) &\mbox{for $x\in (-\infty,u]$,}
\\
\dis\big(1\;\;0\big)M(x)A_1(y) &\mbox{for $x\in [u,v]$,}
\\
\dis\big(1\;\;0\big)N(x)B_{11}(y) &\mbox{for $x\in [v,y]$,}
\\
\dis\big(1\;\;0\big) N(x)B_{12}(y) &\mbox{for $x\in [y,+\infty)$.}
\end{cases}
$$
The matrices are given by
$$
A_1(y)=\g_0(y)M(v)^{-1}N(v) \!\begin{pmatrix}1\\0\end{pmatrix}\!,\quad
B_0(y)=\g_0(y)\!\begin{pmatrix}1\\0\end{pmatrix}\!,
$$
$$
B_{11}(y)=\g_0(y)P_1\!\begin{pmatrix}1\\0\end{pmatrix}\!,\quad
B_{12}(y)=\g_0(y)P_1\!\begin{pmatrix}1\\0\end{pmatrix}-
\kappa(y)N(y)^{-1}\!\begin{pmatrix}0\\1\end{pmatrix}
$$
with
$$
\g_0(y)=\kappa(y)\,\frac{\big(1\;\;0\big) N(y)^{-1}\!\begin{pmatrix}0\\1\end{pmatrix}}{
\big(1\;\;0\big)P_1 \!\begin{pmatrix}1\\0\end{pmatrix}}.
$$

\subsection{Brownian motion}\label{brownian}

In this example, we take for $(X_t)_{t\ge 0}$ linear Brownian motion rescaled
such that $\mathbb{E}(X_t^2)=2t$, that is $X_t=\sqrt2\,B_t$ where
$(B_t)_{t\ge 0}$ is the standard Brownian motion satisfying $\mathbb{E}(B_t^2)=t$.
This choice which corresponds to $\sigma(x)=\sqrt{2}$ and $\tau(x)=0$
is done for simplifying the forthcoming settings.
We take for $a,b,c,d$ the functions
$$
a(x)=\mathrm{e}^{\sqrt{\lambda+\mu}\,x},\quad b(x)=\mathrm{e}^{-\sqrt{\lambda+\mu}\,x},\quad
c(x)=\mathrm{e}^{\sqrt\lambda\,x},\quad d(x)=\mathrm{e}^{-\sqrt\lambda\,x}.
$$
The potential $\rho_{_{\hspace{-.05em}\lambda}}$ writes
$
\rho_{_{\hspace{-.05em}\lambda}}\!(x,y)=\frac{1}{2\sqrt\lambda}\,\mathrm{e}^{-\sqrt\lambda\,|x-y|}
$
and then $w(x)=2\sqrt\lambda$ and $\kappa(x)=1$.
The related matrices $M$ and $N$ write
\begin{align*}
M(x)
&
=\begin{pmatrix} \mathrm{e}^{\sqrt{\lambda+\mu}\,x} & \mathrm{e}^{-\sqrt{\lambda+\mu}\,x} \\
\sqrt{\lambda+\mu}\,\mathrm{e}^{\sqrt{\lambda+\mu}\,x} & -\sqrt{\lambda+\mu}
\,\mathrm{e}^{-\sqrt{\lambda+\mu}\,x} \end{pmatrix}\!,
\\
N(x)
&
=\begin{pmatrix} \mathrm{e}^{\sqrt\lambda\,x} & \mathrm{e}^{-\sqrt\lambda\,x} \\
\sqrt\lambda\,\mathrm{e}^{\sqrt\lambda\,x} & -\sqrt\lambda\,\mathrm{e}^{-\sqrt\lambda\,x} \end{pmatrix}\!,
\end{align*}
and their inverse are given by
\begin{align*}
M(x)^{-1}
&
=\frac{1}{2\sqrt{\lambda+\mu}}\begin{pmatrix} \sqrt{\lambda+\mu}\,\mathrm{e}^{-\sqrt{\lambda+\mu}\,x}
& \mathrm{e}^{-\sqrt{\lambda+\mu}\,x} \\
\sqrt{\lambda+\mu}\,\mathrm{e}^{\sqrt{\lambda+\mu}\,x} & -\mathrm{e}^{\sqrt{\lambda+\mu}\,x} \end{pmatrix}\!,
\\
N(x)^{-1}
&
=\frac{1}{2\sqrt{\lambda}}\begin{pmatrix} \sqrt\lambda\,\mathrm{e}^{-\sqrt\lambda\,x} & \mathrm{e}^{-\sqrt\lambda\,x} \\
\sqrt\lambda\,\mathrm{e}^{\sqrt\lambda\,x} & -\mathrm{e}^{\sqrt\lambda\,x} \end{pmatrix}\!.
\end{align*}
We have
\begin{align*}
M(y)M(x)^{-1}
&
=\begin{pmatrix} \cosh(\sqrt{\lambda+\mu}\,(y-x))& \dis\frac{\sinh(\sqrt{\lambda+\mu}\,(y-x))}{\sqrt{\lambda+\mu}} \\[2ex]
\sqrt{\lambda+\mu}\,\sinh(\sqrt{\lambda+\mu}\,(y-x)) & \cosh(\sqrt{\lambda+\mu}\,(y-x)) \end{pmatrix}
\\
&
=\begin{pmatrix} \cosh(z)& \dis\frac{\sinh(z)}{\sqrt{\lambda+\mu}} \\[2ex]
\sqrt{\lambda+\mu}\,\sinh(z) & \cosh(z) \end{pmatrix}\!,
\end{align*}
where we set $z=\sqrt{\lambda+\mu}\,(y-x)$ for lightening the text.
Next,
\begin{align*}
\lqn{N(y)^{-1}M(y)M(x)^{-1}}
\\[-3ex]
&
=\frac12\begin{pmatrix}
\mathrm{e}^{-\sqrt{\lambda}\,y} \left[\vphantom{\frac{a^2}{a^2}}\right.\!\!\cosh(z)+
\sqrt{\frac{\lambda+\mu}{\lambda}}\sinh(z)\left.\vphantom{\frac{a^2}{a^2}}\!\!\right]
&
\frac{\mathrm{e}^{-\sqrt{\lambda}\,y}}{\sqrt{\lambda}} \left[\vphantom{\frac{a^2}{a^2}}\right.\!\!\cosh(z)+
\sqrt{\frac{\lambda}{\lambda+\mu}}\sinh(z)\left.\vphantom{\frac{a^2}{a^2}}\!\!\right]
\\[2ex]
\mathrm{e}^{\sqrt{\lambda}\,y}\left[\vphantom{\frac{a^2}{a^2}}\right.\!\!\cosh(z)-
\sqrt{\frac{\lambda+\mu}{\lambda}}\sinh(z)\left.\vphantom{\frac{a^2}{a^2}}\!\!\right]
&
\frac{\mathrm{e}^{\sqrt{\lambda}\,y}}{\sqrt{\lambda}}\left[\vphantom{\frac{a^2}{a^2}}\right.\!\!-\cosh(z)+
\sqrt{\frac{\lambda}{\lambda+\mu}}\sinh(z)\left.\vphantom{\frac{a^2}{a^2}}\!\!\right]
\end{pmatrix}\!.
\end{align*}
Making use of the elementary identity
\begin{align*}
\cosh x+a \sinh x=1+2\sinh^2\frac x2+2a\cosh\frac x2\sinh\frac x2
=1+2\sinh\frac x2\left(a\cosh\frac x2+\sinh\frac x2\right)\!,
\end{align*}
we get
\begin{align*}
\lqn{N(y)^{-1}M(y)M(x)^{-1}-N(y)^{-1}}
=\textstyle\sinh(\frac z2)\begin{pmatrix}
\mathrm{e}^{-\sqrt{\lambda}\,y} \left[\vphantom{\frac{a^2}{a^2}}\right.\!\!
\sqrt{\frac{\lambda+\mu}{\lambda}}\cosh(\frac z2)+
\sinh(\frac z2)\left.\vphantom{\frac{a^2}{a^2}}\!\!\right]
&
\frac{\mathrm{e}^{-\sqrt{\lambda}\,y}}{\sqrt{\lambda}} \left[\vphantom{\frac{a^2}{a^2}}\right.\!\!
\sqrt{\frac{\lambda}{\lambda+\mu}}\cosh(\frac z2)+
\sinh(\frac z2)\left.\vphantom{\frac{a^2}{a^2}}\!\!\right]
\\[2ex]
\mathrm{e}^{\sqrt{\lambda}\,y}\left[\vphantom{\frac{a^2}{a^2}}\right.\!\!
-\sqrt{\frac{\lambda+\mu}{\lambda}}\cosh(\frac z2)+
\sinh(\frac z2)\left.\vphantom{\frac{a^2}{a^2}}\!\!\right]
&
\frac{\mathrm{e}^{\sqrt{\lambda}\,y}}{\sqrt{\lambda}} \left[\vphantom{\frac{a^2}{a^2}}\right.\!\!
\sqrt{\frac{\lambda}{\lambda+\mu}}\cosh(\frac z2)-
\sinh(\frac z2)\left.\vphantom{\frac{a^2}{a^2}}\!\!\right]
\end{pmatrix}\!.
\end{align*}
Now,
\begin{align*}
\lqn{N(y)^{-1}M(y)M(x)^{-1}N(x)}
\\[-3ex]
&\hspace*{-5em}
=\frac12\left(\begin{matrix} \mathrm{e}^{\sqrt{\lambda}\,(x-y)}
\left[2\cosh(z)+\left(\vphantom{\frac{a^2}{a^2}}\right.\right.\!\!\!
\sqrt{\frac{\lambda+\mu}{\lambda}}+\sqrt{\frac{\lambda}{\lambda+\mu}}\!\left.\left.\vphantom{\frac{a^2}{a^2}}\right)
\sinh(z)\right] \\[2ex]
\mathrm{e}^{\sqrt{\lambda}\,(x+y)}\left(\vphantom{\frac{a^2}{a^2}}\right.\!\!
\sqrt{\frac{\lambda}{\lambda+\mu}}-\sqrt{\frac{\lambda+\mu}{\lambda}}\left.\vphantom{\frac{a^2}{a^2}}\right)\sinh(z)
\end{matrix}\right.
\\[3ex]
&
\hspace{10em}\left.\begin{matrix}
\mathrm{e}^{-\sqrt{\lambda}\,(x+y)}\left(\vphantom{\frac{a^2}{a^2}}\right.\!\!
\sqrt{\frac{\lambda+\mu}{\lambda}}-\sqrt{\frac{\lambda}{\lambda+\mu}}\left.\vphantom{\frac{a^2}{a^2}}\right)\sinh(z)
\\[2ex]
\mathrm{e}^{\sqrt{\lambda}\,(y-x)}\left[2\cosh(z)-\left(\vphantom{\frac{a^2}{a^2}}\right.\right.\!\!\!
\sqrt{\frac{\lambda+\mu}{\lambda}}+\sqrt{\frac{\lambda}{\lambda+\mu}}\!\left.\left.\vphantom{\frac{a^2}{a^2}}\right)\sinh(z)\right]
\end{matrix}\right)\!.
\end{align*}
Observing that
\begin{align*}
2\cosh x+\left(a+\frac 1a\right)\sinh x
&
=2\left(\cosh^2\frac x2+\sinh^2\frac x2\right)
+2\left(a+\frac 1a\right)\cosh \frac x2\sinh \frac x2
\\
&
=\frac 2a\left(a\cosh \frac x2+\sinh \frac x2\right) \!\!\left(\cosh \frac x2+a\sinh \frac x2\right)\!,
\end{align*}
we obtain
\begin{align*}
\lqn{
N(y)^{-1}M(y)M(x)^{-1}N(x)
}
\\[-3ex]
&\hspace*{-5em}
=\frac{1}{\sqrt{\lambda(\lambda+\mu)}}\left(\begin{matrix} \mathrm{e}^{\sqrt{\lambda}\,(x-y)}
\big[\sqrt{\lambda}\cosh(\frac z2)+\sqrt{\lambda+\mu}\sinh(\frac z2)\big]\!
\big[\sqrt{\lambda+\mu}\cosh(\frac z2)+\sqrt{\lambda}\sinh(\frac z2)\big]
\\[2ex]
-\mu\mathrm{e}^{\sqrt{\lambda}\,(x+y)} \cosh(\frac z2)\sinh(\frac z2)
\end{matrix}\right.
\\[3ex]
&
\hspace{2em}\left.\begin{matrix}
\mu\mathrm{e}^{-\sqrt{\lambda}\,(x+y)}\cosh(\frac z2)\sinh(\frac z2)
\\[2ex]
\mathrm{e}^{\sqrt{\lambda}\,(y-x)}
\big[\sqrt{\lambda}\cosh(\frac z2)-\sqrt{\lambda+\mu}\sinh(\frac z2)\big]\!
\big[\sqrt{\lambda+\mu}\cosh(\frac z2)-\sqrt{\lambda}\sinh(\frac z2)\big]
\end{matrix}\right)\!.
\end{align*}

\subsection{Case of one bounded interval for Brownian motion}

We now consider the sojourn time of Brownian motion in the interval $[u,v]$.
In order to lighten the paper, we only compute the distribution of $T_t$,
that of $(T_t,X_t)$ being more cumbersome.

For evaluating $\g_0$, we compute, with $w=\frac12\sqrt{\lambda+\mu}\,(v-u)$,
$$
\big(1\;\;0\big)P_1\!\begin{pmatrix}1\\0\end{pmatrix}
=\frac{\mathrm{e}^{\sqrt{\lambda}\,(u-v)}}{\sqrt{\lambda(\lambda+\mu)}}
\big[\sqrt{\lambda}\cosh(w)+\sqrt{\lambda+\mu}\sinh(w)\big]\!
\big[\sqrt{\lambda+\mu}\cosh(w)+\sqrt{\lambda}\sinh(w)\big]
$$
and
$$
\big(1\;\;0\big)Q_1\!\begin{pmatrix}1\\0\end{pmatrix}
=\frac{\mathrm{e}^{-\sqrt{\lambda}\,v}}{\sqrt{\lambda}}
\sinh(w) \big[\sqrt{\lambda+\mu}\cosh(w)+\sqrt{\lambda}\sinh(w)\big]
$$
from which we deduce
$$
\g_0=-\frac{\mu}{\lambda(\lambda+\mu)}\,\frac{\big(1\;\;0\big)Q_1\!\begin{pmatrix}1\\0\end{pmatrix}}{
\big(1\;\;0\big)P_1\!\begin{pmatrix}1\\0\end{pmatrix}}
=-\frac{\mu}{\lambda\sqrt{\lambda+\mu}}\,\frac{\mathrm{e}^{-\sqrt{\lambda}\,u}\sinh(w)}{
\sqrt{\lambda}\cosh(w)+\sqrt{\lambda+\mu}\sinh(w)}.
$$
For evaluating $\d_1$, we compute
$$
\begin{pmatrix}0\;\;1\end{pmatrix}\!\tilde{P}_1\!\begin{pmatrix}0\\1\end{pmatrix}
=\frac{\mathrm{e}^{\sqrt{\lambda}\,(u-v)}}{\sqrt{\lambda(\lambda+\mu)}}
\big[\sqrt{\lambda}\cosh(w)+\sqrt{\lambda+\mu}\sinh(w)\big]\!
\big[\sqrt{\lambda+\mu}\cosh(w)+\sqrt{\lambda}\sinh(w)\big]
$$
and
$$
\begin{pmatrix}0\;\;1\end{pmatrix}\!\tilde{Q}_1\!\begin{pmatrix}1\\0\end{pmatrix}
=\frac{\mathrm{e}^{\sqrt{\lambda}\,u}}{\sqrt{\lambda}}
\sinh(w) \big[\sqrt{\lambda+\mu}\cosh(w)+\sqrt{\lambda}\sinh(w)\big]
$$
from which we deduce
$$
\d_1=-\frac{\mu}{\lambda(\lambda+\mu)}\,\frac{\begin{pmatrix}0\;\;1\end{pmatrix}
\!\tilde{Q}_1\!\begin{pmatrix}1\\0\end{pmatrix}}{
\begin{pmatrix}0\;\;1\end{pmatrix}\!\tilde{P}_1\!\begin{pmatrix}0\\1\end{pmatrix}}
=-\frac{\mu}{\lambda\sqrt{\lambda+\mu}}\,\frac{\mathrm{e}^{\sqrt{\lambda}\,v}\sinh(w)}{
\sqrt{\lambda}\cosh(w)+\sqrt{\lambda+\mu}\sinh(w)}.
$$
This yields
$$
\big(1\;\;0\big) N(x)B_0 = \g_0\big(1\;\;0\big) N(x)\!\begin{pmatrix}1\\0\end{pmatrix}
=-\frac{\mu}{\lambda\sqrt{\lambda+\mu}}\,\frac{\sinh(w)\mathrm{e}^{\sqrt{\lambda}\,(x-u)}}{
\sqrt{\lambda}\cosh(w)+\sqrt{\lambda+\mu}\sinh(w)}
$$
$$
\big(1\;\;0\big) N(x)B_1 = \d_1\big(1\;\;0\big) N(x)\!\begin{pmatrix}0\\1\end{pmatrix}
=-\frac{\mu}{\lambda\sqrt{\lambda+\mu}}\,\frac{\sinh(w)\mathrm{e}^{\sqrt{\lambda}\,(v-x)}}{
\sqrt{\lambda}\cosh(w)+\sqrt{\lambda+\mu}\sinh(w)}.
$$

On the other hand,
$$
M(u)^{-1}\!\begin{pmatrix}1\\0\end{pmatrix}
=\frac12\,\mathrm{e}^{\sqrt{\lambda}\,u}
\begin{pmatrix} \mathrm{e}^{-\sqrt{\lambda+\mu}\,u}
\\
\mathrm{e}^{\sqrt{\lambda+\mu}\,u} \end{pmatrix}\!,
$$
\begin{align*}
M(u)^{-1}N(u)\!\begin{pmatrix}1\\0\end{pmatrix}
&
=\frac{\mathrm{e}^{\sqrt{\lambda}\,u} }{2\sqrt{\lambda+\mu}}
\begin{pmatrix} \sqrt{\lambda+\mu}\,\mathrm{e}^{-\sqrt{\lambda+\mu}\,u}
& \mathrm{e}^{-\sqrt{\lambda+\mu}\,u} \\
\sqrt{\lambda+\mu}\,\mathrm{e}^{\sqrt{\lambda+\mu}\,u}
& -\mathrm{e}^{\sqrt{\lambda+\mu}\,u} \end{pmatrix}
\!\!\begin{pmatrix}1\\ \sqrt{\lambda}\end{pmatrix}
\\
&
=\frac{\mathrm{e}^{\sqrt{\lambda}\,u} }{2\sqrt{\lambda+\mu}}
\begin{pmatrix} \big[\sqrt{\lambda+\mu}+\sqrt{\lambda}\,\big]\,\mathrm{e}^{-\sqrt{\lambda+\mu}\,u}
\\[.5ex]
\big[\sqrt{\lambda+\mu}-\sqrt{\lambda}\,\big]\,\mathrm{e}^{\sqrt{\lambda+\mu}\,u} \end{pmatrix}
\end{align*}
and then
\begin{align*}
A_1
&
=\g_0 M(u)^{-1}N(u)\!\begin{pmatrix}1\\0\end{pmatrix}
+\frac{\mu}{\lambda(\lambda+\mu)}M(u)^{-1}\!\begin{pmatrix}1\\0\end{pmatrix}
\\
&
=\frac{\mu}{2\lambda(\lambda+\mu)}
\begin{pmatrix} \dis\bigg[1-\frac{\big(\sqrt{\lambda+\mu}+\sqrt{\lambda}\,\big)\sinh(w)}{
\sqrt{\lambda}\cosh(w)+\sqrt{\lambda+\mu}\sinh(w)}\bigg]\mathrm{e}^{-\sqrt{\lambda+\mu}\,u}
\\[.5ex]
\dis\bigg[1-\frac{\big(\sqrt{\lambda+\mu}-\sqrt{\lambda}\,\big)\sinh(w)}{
\sqrt{\lambda}\cosh(w)+\sqrt{\lambda+\mu}\sinh(w)}\bigg]\mathrm{e}^{\sqrt{\lambda+\mu}\,u}
\end{pmatrix}
\\
&
=\frac{\mu}{2\sqrt{\lambda}\,(\lambda+\mu)}\,\frac{1}{\sqrt{\lambda}\cosh(w)+\sqrt{\lambda+\mu}\sinh(w)}
\begin{pmatrix} [\cosh(w)-\sinh(w)]\,\mathrm{e}^{-\sqrt{\lambda+\mu}\,u}
\\
[\cosh(w)+\sinh(w)]\,\mathrm{e}^{\sqrt{\lambda+\mu}\,u}
\end{pmatrix}
\\
&
=\frac{\mu}{2\sqrt{\lambda}\,(\lambda+\mu)}\,\frac{1}{\sqrt{\lambda}\cosh(w)+\sqrt{\lambda+\mu}\sinh(w)}
\begin{pmatrix} \mathrm{e}^{-\sqrt{\lambda+\mu}\,(u+v)/2}
\\
\mathrm{e}^{\sqrt{\lambda+\mu}\,(u+v)/2}
\end{pmatrix}\!.
\end{align*}
Now,
\begin{align*}
\big(1\;\;0\big) M(x)A_1
&
=\frac{\mu}{2\sqrt{\lambda}\,(\lambda+\mu)}\,\frac{1}{\sqrt{\lambda}\cosh(w)+\sqrt{\lambda+\mu}\sinh(w)}
\begin{pmatrix} \mathrm{e}^{\sqrt{\lambda+\mu}\,x} & \mathrm{e}^{-\sqrt{\lambda+\mu}\,x} \end{pmatrix}
\!\begin{pmatrix} \mathrm{e}^{-\sqrt{\lambda+\mu}\,(u+v)/2}
\\
\mathrm{e}^{\sqrt{\lambda+\mu}\,(u+v)/2}
\end{pmatrix}
\\
&
=\frac{\mu}{\sqrt{\lambda}\,(\lambda+\mu)}\,\frac{\cosh\!\big(\sqrt{\lambda+\mu}\,(x-(u+v)/2)\big)}{
\sqrt{\lambda}\cosh(w)+\sqrt{\lambda+\mu}\sinh(w)}.
\end{align*}

We sum up below the results we have obtained.
%
\begin{theo}
The distribution of the sojourn time $T_t$ inside the interval $[u,v]$ for
Brownian motion is characterized by
$$
\int_0^{\infty} \mathrm{e}^{-\lambda t}\, \mathbb{E}_x\!\left(\mathrm{e}^{-\mu T_t}\right) \mathrm{d} t
=\begin{cases}
\dis\frac{1}{\lambda}\bigg[1-\frac{\mu}{\sqrt{\lambda+\mu}}\,\frac{\sinh(w)\mathrm{e}^{\sqrt{\lambda}\,(x-u)}}{
\sqrt{\lambda}\cosh(w)+\sqrt{\lambda+\mu}\sinh(w)}\bigg]
&\mbox{for $x\in (-\infty,u]$,}
\\[2ex]
\dis
\dis\frac{1}{\lambda+\mu}\bigg[1-\frac{\mu}{\sqrt{\lambda}}\,\frac{\cosh\!\big(\sqrt{\lambda+\mu}\,(x-(u+v)/2)\big)}{
\sqrt{\lambda}\cosh(w)+\sqrt{\lambda+\mu}\sinh(w)}\bigg]
&\mbox{for $x\in [u,v]$,}
\\[2ex]
\dis\frac{1}{\lambda}\bigg[1-\frac{\mu}{\sqrt{\lambda+\mu}}\,\frac{\sinh(w)\mathrm{e}^{\sqrt{\lambda}\,(v-x)}}{
\sqrt{\lambda}\cosh(w)+\sqrt{\lambda+\mu}\sinh(w)}\bigg]
&\mbox{for $x\in [v,+\infty)$,}
\end{cases}
$$
where $w=\frac12\sqrt{\lambda+\mu}\,(v-u)$.
\end{theo}
%
We retrieve the well-known distribution (1.7.1), p. 140 of~\cite{borodin}.

\subsection{Local time in a finite set for Brownian motion}

Let us apply our results to the following set ($u_1,\dots,u_n$ are real numbers
such that $u_1<\dots<u_n$):
$$
E_{\e}=\bigcup_{i=1}^n [u_i-\e,u_i+\e]
$$
where $\e>0$ is subject to tend to $0$ and denote
$T_{t,\e}=\int_0^t \ind_{E_{\e}}(X_s)\,\mathrm{d}s$.
Set also $u_0=-\infty$ and $u_{n+1}=+\infty$.
The local time in the set $\{u_1,\dots,u_n\}$ of Brownian motion $(X_t)_{t\ge 0}$
up to time~$t$ is defined by
$$
L_t=\lim_{\e\to 0}\frac{1}{\e}\,T_{t,\e}.
$$
Of course, as previously, we can decompose $L_t$ into the sum
$$
L_t=\sum_{i=1}^nL_t^i
$$
where, for any $i\in\{1,\dots,n\}$,
$$
L_t^i=\lim_{\e\to 0}\frac{1}{\e}\int_0^t \ind_{[u_i-\e,u_i+\e]}(X_s)\,\mathrm{d}s.
$$
Let us introduce the vector of local times at each point $u_i$:
$\mathbf{L}_t=(L_t^1,\dots,L_t^n)$.
Set, for $i\in\{1,\dots,n\}$,
\begin{align*}
\nu_{i,\e}
&
=\frac{\mu_i/\e}{\lambda(\lambda+\mu_i/\e)},
\\
P_{i,\e}
&
=N(u_i+\e)^{-1}M_i(u_i+\e)M_i(u_i-\e)^{-1}N(u_i-\e),
\\
Q_{i,\e}
&
=N(u_i+\e)^{-1}M_i(u_i+\e)M_i(u_i-\e)^{-1}-N(u_i+\e)^{-1}.
\\
R_{i,\e}
&
=P_{i,\e}P_{i-1,\e}\dots P_{1,\e},
\\
S_{i,\e}
&
=\nu_{i,\e}Q_{i,\e}+\nu_{i-1,\e}P_{i,\e}Q_{i-1,\e}+\nu_{i-2,\e}P_{i,\e}P_{i-1,\e}Q_{i-2,\e}
+\dots +\nu_{1,\e}P_{i,\e}P_{i-1,\e}\dots P_{2,\e}Q_{1,\e},
\\
B_{i,\e}
&
=S_{i,\e} C_0-\frac{\big(1\;\;0\big)S_{n,\e} C_0}{\big(1\;\;0\big)R_{n,\e} C_0}R_{i,\e} C_0.
\end{align*}
As usual, we set $R_{0,\varepsilon}=I$, $S_{0,\varepsilon}=O$ and
$B_{0,\varepsilon}=-\frac{(1\;0)S_{n,\e} C_0}{(1\;0)R_{n,\e} C_0}\,C_0$.

In this part, we compute the following limit, for $x\in (u_i,u_{i+1})$
(which entails that $x\in (u_i+\e,u_{i+1}-\e)$ for small enough $\e$) and $i\in\{0,1,\dots,n\}$,
$$
\int_0^{\infty} \mathrm{e}^{-\lambda t}\, \mathbb{E}_x\!
\left(\mathrm{e}^{-<\boldsymbol{\mu},\mathbf{L}_t>}\right) \mathrm{d} t
=\lim_{\e\to 0}\big(1\;\;0\big) N(x)B_{i,\e}+\frac{1}{\lambda}
$$
where
$$
N(x)=\begin{pmatrix} \mathrm{e}^{\sqrt\lambda\,x} & \mathrm{e}^{-\sqrt\lambda\,x} \\
\sqrt\lambda\,\mathrm{e}^{\sqrt\lambda\,x} & -\sqrt\lambda\,\mathrm{e}^{-\sqrt\lambda\,x} \end{pmatrix}\!.
$$
Set $\epsilon=\e\sqrt{\lambda+\mu_i/\e}$. Using the results of
Section~\ref{brownian}, we have for $i\in\{1,\dots,n\}$
\begin{align*}
P_{i,\e}
&=\frac{1}{\sqrt{\lambda(\lambda+\mu_i/\e)}}
\left(\begin{matrix} \mathrm{e}^{-2\e\sqrt{\lambda}}
\big[\sqrt{\lambda}\cosh(\epsilon)+\sqrt{\lambda+\mu_i/\e}\sinh(\epsilon)\big]\!
\big[\sqrt{\lambda+\mu_i/\e}\cosh(\epsilon)+\sqrt{\lambda}\sinh(\epsilon)\big]
\\[2ex]
-(\mu_i/\e)\mathrm{e}^{2u_i\sqrt{\lambda}} \cosh(\epsilon)\sinh(\epsilon)
\end{matrix}\right.
\\[3ex]
&
\hspace{10em}\left.\begin{matrix}
(\mu_i/\e)\mathrm{e}^{-2u_i\sqrt{\lambda}}\cosh(\epsilon)\sinh(\epsilon)
\\[2ex]
\mathrm{e}^{2\e\sqrt{\lambda}}
\big[\sqrt{\lambda}\cosh(\epsilon)-\sqrt{\lambda+\mu_i/\e}\sinh(\epsilon)\big]\!
\big[\sqrt{\lambda+\mu_i/\e}\cosh(\epsilon)-\sqrt{\lambda}\sinh(\epsilon)\big]
\end{matrix}\right)
\end{align*}
and
\begin{align*}
Q_{i,\e}
=\sinh(\epsilon)\begin{pmatrix}
\mathrm{e}^{-\sqrt{\lambda}\,(u_i+\e)} \left[\vphantom{\frac{a^2}{a^2}}\right.\!\!
\sqrt{\frac{\lambda+\mu_i/\e}{\lambda}}\cosh(\epsilon)+
\sinh(\epsilon)\left.\vphantom{\frac{a^2}{a^2}}\!\!\right]
&
\frac{\mathrm{e}^{-\sqrt{\lambda}\,(u_i+\e)}}{\sqrt{\lambda}}
\left[\vphantom{\frac{a^2}{a^2}}\right.\!\!
\sqrt{\frac{\lambda}{\lambda+\mu_i/\e}}\cosh(\epsilon)+
\sinh(\epsilon)\left.\vphantom{\frac{a^2}{a^2}}\!\!\right]
\\[2ex]
\mathrm{e}^{\sqrt{\lambda}\,(u_i+\e)}\left[\vphantom{\frac{a^2}{a^2}}\right.\!\!
-\sqrt{\frac{\lambda+\mu_i/\e}{\lambda}}\cosh(\epsilon)+
\sinh(\epsilon)\left.\vphantom{\frac{a^2}{a^2}}\!\!\right]
&
\frac{\mathrm{e}^{\sqrt{\lambda}\,(u_i+\e)}}{\sqrt{\lambda}} \left[\vphantom{\frac{a^2}{a^2}}\right.\!\!
\sqrt{\frac{\lambda}{\lambda+\mu_i/\e}}\cosh(\epsilon)-
\sinh(\epsilon)\left.\vphantom{\frac{a^2}{a^2}}\!\!\right]
\end{pmatrix}\!.
\end{align*}
By the elementary asymptotics for $\e\to0^+$
\begin{align*}
\sqrt{\lambda}\cosh(\epsilon)+\sqrt{\lambda+\mu_i/\e}\sinh(\epsilon)
&
\sim\sqrt{\lambda}+\mu_i,
\\
\sqrt{\lambda+\mu_i/\e}\cosh(\epsilon)+\sqrt{\lambda}\sinh(\epsilon)
&
\sim\sqrt{\mu_i/\e},
\\
\cosh(\epsilon)\sinh(\epsilon)
&
\sim\sqrt{\mu_i\e},
\end{align*}
we get, for $i\in\{1,\dots,n\}$, $\lim_{\e\to0^+} P_{i,\e}=\bar{P}_i$,
$\lim_{\e\to0^+} Q_{i,\e}=\bar{Q}_i$, $\lim_{\e\to0^+} R_{i,\e}=\bar{R}_i$,
$\lim_{\e\to0^+} S_{i,\e}=\bar{S}_i$ with
\begin{equation}\label{Pbar}
\bar{P}_i=\frac{1}{\sqrt{\lambda}}\begin{pmatrix}
\sqrt{\lambda}+\mu_i & \mu_i\mathrm{e}^{-2\sqrt{\lambda}\,u_i}
\\
-\mu_i\mathrm{e}^{2\sqrt{\lambda}\,u_i} & \sqrt{\lambda}-\mu_i
\end{pmatrix}\!,\quad
\bar{Q}_i=\frac{\mu_i}{\sqrt{\lambda}}
\begin{pmatrix}
\mathrm{e}^{-\sqrt{\lambda}\,u_i}&0
\\
-\mathrm{e}^{\sqrt{\lambda}\,u_i}&0
\end{pmatrix}\!,
\end{equation}
\begin{equation}\label{Rbar}
\bar{R}_i=\bar{P}_i\bar{P}_{i-1}\dots \bar{P}_1,\quad
\bar{S}_i=\frac{1}{\lambda}\left[\bar{Q}_i+\bar{P}_i\bar{Q}_{i-1}+\bar{P}_i\bar{P}_{i-1}
\bar{Q}_{i-2}+\dots +\bar{P}_i\bar{P}_{i-1}\dots \bar{P}_2\bar{Q}_1\right]\!.
\end{equation}
We also have, for $i\in\{0,1,\dots,n\}$, $\lim_{\e\to0^+} B_{i,\e}=\bar{B}_i$ where
\begin{equation}\label{Bbar}
\bar{B}_i=\bar{S}_i\!\begin{pmatrix}1\\0\end{pmatrix}
-\bar{\gamma}_0\bar{R}_i\!\begin{pmatrix}1\\0\end{pmatrix}
\quad\mbox{with}\quad
\bar{\gamma}_0=\frac{\big(1\;\;0\big)\bar{S}_n\!\begin{pmatrix}1\\ 0\end{pmatrix}}{
\big(1\;\;0\big)\bar{R}_n\!\begin{pmatrix}1\\ 0\end{pmatrix}}.
\end{equation}
We can now state the following result.
%
\begin{theo}
The iterated Laplace transform of the vector of Brownian
local times $\mathbf{L}_t$ at points $u_1,\dots,u_n$ is given, for $x\in (u_i,u_{i+1})$
and $i\in\{0,1,\dots,n\}$, by
$$
\int_0^{\infty} \mathrm{e}^{-\lambda t}\, \mathbb{E}_x\!
\left(\mathrm{e}^{-<\boldsymbol{\mu},\mathbf{L}_t>}\right) \mathrm{d} t
=\big(1\;\;0\big) N(x)\bar{B}_i+\frac{1}{\lambda}
$$
where the matrices $\bar{B}_i$, $i\in\{0,1,\dots,n\}$, are given by~(\ref{Pbar}),
(\ref{Rbar}) and~(\ref{Bbar}).
Additionally, this formula holds also for $x=u_i$, $i\in\{1,\dots,n\}$.
In particular, the iterated Laplace transform of the Brownian local time $L_t$ in
$\{u_1,\dots,u_n\}$, namely $\int_0^{\infty} \mathrm{e}^{-\lambda t}\, \mathbb{E}_x\!
\left(\mathrm{e}^{-\mu L_t}\right) \mathrm{d} t$, can be deduced from the
previous formula by choosing $\boldsymbol{\mu}=(\mu,\dots,\mu)$.
\end{theo}
%
\dem
It remains to prove the assertion concerning the case $x=u_i$.
Observing that we have
$$
\big(1\;\;0\big) N(u_i)\bar{P}_i
=\frac{1}{\sqrt{\lambda}}\begin{pmatrix}
\mathrm{e}^{\sqrt{\lambda}\,u_i}&\mathrm{e}^{-\sqrt{\lambda}\,u_i}
\end{pmatrix}\!\!
\begin{pmatrix}
\sqrt{\lambda}+\mu_i & \mu_i\mathrm{e}^{-2\sqrt{\lambda}\,u_i}
\\
-\mu_i\mathrm{e}^{2\sqrt{\lambda}\,u_i} & \sqrt{\lambda}-\mu_i
\end{pmatrix}
=\begin{pmatrix}
\mathrm{e}^{\sqrt{\lambda}\,u_i}&\mathrm{e}^{-\sqrt{\lambda}\,u_i}
\end{pmatrix}
=\big(1\;\;0\big) N(u_i),
$$
and
$$
\big(1\;\;0\big) N(u_i)\bar{Q}_i
=\frac{\mu_i}{\sqrt{\lambda}}
\begin{pmatrix}
\mathrm{e}^{\sqrt{\lambda}\,u_i}&\mathrm{e}^{-\sqrt{\lambda}\,u_i}
\end{pmatrix}\!\!
\begin{pmatrix}
\mathrm{e}^{-\sqrt{\lambda}\,u_i}&0
\\
-\mathrm{e}^{\sqrt{\lambda}\,u_i}&0
\end{pmatrix}
=\big(0\;\;0\big),
$$
we deduce that
$$
\big(1\;\;0\big) N(u_i)\bar{R}_i=\big(1\;\;0\big) N(u_i)\bar{P}_i\bar{P}_{i-1}\dots \bar{P}_1
=\big(1\;\;0\big) N(u_i)\bar{P}_{i-1}\dots \bar{P}_1=\big(1\;\;0\big) N(u_i)\bar{R}_{i-1}
$$
and
\begin{align*}
\big(1\;\;0\big) N(u_i)\bar{S}_i
&
=\frac{1}{\lambda}\left[\big(1\;\;0\big) N(u_i)\bar{Q}_i+\big(1\;\;0\big) N(u_i)\bar{P}_i(\bar{Q}_{i-1}+\dots+
\bar{P}_{i-1}\dots \bar{P}_2\bar{Q}_1)\right]
\\
&
=\frac{1}{\lambda}\big(1\;\;0\big) N(u_i)(\bar{Q}_{i-1}+\dots+\bar{P}_{i-1}\dots \bar{P}_2\bar{Q}_1)
=\big(1\;\;0\big) N(u_i)\bar{S}_{i-1}.
\end{align*}
As a result, since $\bar{B}_i$ is linear combination of the matrices
$\bar{R}_i$ and $\bar{S}_i$, we have
$$
\big(1\;\;0\big) N(u_i)\bar{B}_i=\big(1\;\;0\big) N(u_i)\bar{B}_{i-1}.
$$
This proves that
$$
\lim_{x\to u_i^+}\int_0^{\infty} \mathrm{e}^{-\lambda t}\, \mathbb{E}_x
\!\left(\mathrm{e}^{-<\boldsymbol{\mu},\mathbf{L}_t>}\right) \mathrm{d} t
=\lim_{x\to u_i^-}\int_0^{\infty} \mathrm{e}^{-\lambda t}\, \mathbb{E}_x
\!\left(\mathrm{e}^{-<\boldsymbol{\mu},\mathbf{L}_t>}\right) \mathrm{d} t
$$
and then, by continuity with respect to $x$,
$$
\int_0^{\infty} \mathrm{e}^{-\lambda t}\, \mathbb{E}_{u_i}
\!\left(\mathrm{e}^{-<\boldsymbol{\mu},\mathbf{L}_t>}\right) \mathrm{d} t
=\lim_{x\to u_i}\int_0^{\infty} \mathrm{e}^{-\lambda t}\, \mathbb{E}_x
\!\left(\mathrm{e}^{-<\boldsymbol{\mu},\mathbf{L}_t>}\right) \mathrm{d} t
=\big(1\;\;0\big) N(u_i)\bar{B}_i+\frac{1}{\lambda}.
$$
\fin
%

We end up this part be considering the particular cases $n=1$ and $n=2$.

\vspace{\baselineskip}
\noindent\textsl{\textbf{Local time in $\{u\}$}}
\vspace{\baselineskip}

Suppose that $n=1$ and set $u_1=u$, $\mu_1=\mu$.
We have $\bar{R}_0=I,\bar{S}_0=O,\bar{R}_1=\bar{P}_1,\bar{S}_1=\frac{1}{\lambda}\,\bar{Q}_1$.
Therefore,
$$
\bar{R}_1=\frac{1}{\sqrt{\lambda}}\begin{pmatrix}
\sqrt{\lambda}+\mu & \mu\mathrm{e}^{-2\sqrt{\lambda}\,u}
\\
-\mu\mathrm{e}^{2\sqrt{\lambda}\,u} & \sqrt{\lambda}-\mu
\end{pmatrix}\!,\quad
\bar{S}_1=\frac{\mu}{\lambda^{3/2}}
\begin{pmatrix}
\mathrm{e}^{-\sqrt{\lambda}\,u}&0
\\
-\mathrm{e}^{\sqrt{\lambda}\,u}&0
\end{pmatrix}\!,
$$
$$
\bar{B}_0=\bar{S}_0\!\begin{pmatrix}1\\0\end{pmatrix}
-\frac{\big(1\;\;0\big)\bar{S}_1\!\begin{pmatrix}1\\ 0\end{pmatrix}}{
\big(1\;\;0\big)\bar{R}_1\!\begin{pmatrix}1\\ 0\end{pmatrix}}
\bar{R}_0\!\begin{pmatrix}1\\0\end{pmatrix}
=-\frac{\mu\mathrm{e}^{-\sqrt{\lambda}\,u}}{\lambda(\sqrt\lambda+\mu)}\begin{pmatrix}1\\0\end{pmatrix}
$$
and
$$
\bar{B}_1=\bar{S}_1\!\begin{pmatrix}1\\0\end{pmatrix}
-\frac{\big(1\;\;0\big)\bar{S}_1\!\begin{pmatrix}1\\ 0\end{pmatrix}}{
\big(1\;\;0\big)\bar{R}_1\!\begin{pmatrix}1\\ 0\end{pmatrix}}
\bar{R}_1\!\begin{pmatrix}1\\0\end{pmatrix}
=\frac{\mu}{\lambda^{3/2}}\left[\!
\begin{pmatrix}
\mathrm{e}^{-\sqrt{\lambda}\,u}
\\
-\mathrm{e}^{\sqrt{\lambda}\,u}
\end{pmatrix}
-\frac{\mathrm{e}^{-\sqrt{\lambda}\,u}}{\sqrt\lambda+\mu}
\begin{pmatrix}
\sqrt{\lambda}+\mu
\\
-\mu\mathrm{e}^{2\sqrt{\lambda}\,u}
\end{pmatrix}\!\right]
=-\frac{\mu\mathrm{e}^{\sqrt{\lambda}\,u}}{\lambda(\sqrt\lambda+\mu)}
\begin{pmatrix}0\\1\end{pmatrix}\!.
$$
Next,
$$
\big(1\;\;0\big) N(x)\bar{B}_0=-\frac{\mu}{\lambda(\sqrt\lambda+\mu)}\,
\mathrm{e}^{\sqrt{\lambda}\,(x-u)},\quad
\big(1\;\;0\big) N(x)\bar{B}_1=-\frac{\mu}{\lambda(\sqrt\lambda+\mu)}\,
\mathrm{e}^{\sqrt{\lambda}\,(u-x)}.
$$
As a result, the iterated Laplace transform of the local time $L_t$ in $\{u\}$
is given, for any $x\in \mathbb{R}$, by
$$
\int_0^{\infty} \mathrm{e}^{-\lambda t}\, \mathbb{E}_x\!
\left(\mathrm{e}^{-\mu L_t}\right) \mathrm{d} t
=\frac{1}{\lambda}\bigg[1-\frac{\mu}{\sqrt\lambda+\mu}\,
\mathrm{e}^{-\sqrt{\lambda}\,|x-u|}\bigg].
$$
We retrieve formula~(1.3.1), p. 126 of~\cite{borodin}.

\vspace{\baselineskip}
\noindent\textsl{\textbf{Local time in $\{u,v\}$}}
\vspace{\baselineskip}

Suppose that $n=2$ and set $u_1=u,u_2=v,\mu_1=\mu,\mu_2=\nu$. We have
$\bar{R}_0=I,\bar{S}_0=O,\bar{R}_1=\bar{P}_1,\bar{S}_1=\frac{1}{\lambda}\,\bar{Q}_1,
\bar{R}_2=\bar{P}_2\bar{P}_1,\bar{S}_2=\frac{1}{\lambda}\left(\bar{Q}_2+\bar{P}_2\bar{Q}_1\right)$,
$\bar{\gamma}_0=\big[\big(1\;\;0\big)\bar{S}_2\!\begin{pmatrix}1\\ 0\end{pmatrix}\!\big]
/\big[\big(1\;\;0\big)\bar{R}_2\!\begin{pmatrix}1\\ 0\end{pmatrix}\!\big]$.
Explicitly,
$$
\bar{R}_1=\frac{1}{\sqrt{\lambda}}\begin{pmatrix}
\sqrt{\lambda}+\mu & \mu\mathrm{e}^{-2\sqrt{\lambda}\,u}
\\
-\mu\mathrm{e}^{2\sqrt{\lambda}\,u} & \sqrt{\lambda}-\mu
\end{pmatrix}\!,\quad
\bar{S}_1=\frac{\mu}{\lambda^{3/2}}
\begin{pmatrix}
\mathrm{e}^{-\sqrt{\lambda}\,u}&0
\\
-\mathrm{e}^{\sqrt{\lambda}\,u}&0
\end{pmatrix}\!,
$$
\begin{align*}
\bar{R}_2
&
=\frac{1}{\lambda}\begin{pmatrix}
\sqrt{\lambda}+\nu & \nu\mathrm{e}^{-2\sqrt{\lambda}\,v}
\\
-\nu\mathrm{e}^{2\sqrt{\lambda}\,v} & \sqrt{\lambda}-\nu
\end{pmatrix}\!\!
\begin{pmatrix}
\sqrt{\lambda}+\mu & \mu\mathrm{e}^{-2\sqrt{\lambda}\,u}
\\
-\mu\mathrm{e}^{2\sqrt{\lambda}\,u} & \sqrt{\lambda}-\mu
\end{pmatrix}
\\
&
=\frac{1}{\lambda}\begin{pmatrix}
\big(\sqrt{\lambda}+\mu\big)\!\big(\sqrt{\lambda}+\nu\big)
-\mu\nu \mathrm{e}^{2\sqrt{\lambda}\,(u-v)}
&
\big(\sqrt{\lambda}+\nu\big)\mu\mathrm{e}^{-2\sqrt{\lambda}\,u}
+\big(\sqrt{\lambda}-\mu\big)\nu\mathrm{e}^{-2\sqrt{\lambda}\,v}
\\
-\big(\sqrt{\lambda}+\nu\big)\mu\mathrm{e}^{2\sqrt{\lambda}\,u}
-\big(\sqrt{\lambda}+\mu\big)\nu\mathrm{e}^{2\sqrt{\lambda}\,v}
&
\big(\sqrt{\lambda}-\mu\big)\!\big(\sqrt{\lambda}-\nu\big)
-\mu\nu \mathrm{e}^{2\sqrt{\lambda}\,(v-u)}
\end{pmatrix}\!,
\end{align*}
\begin{align*}
\bar{S}_2
&
=\frac{\nu}{\lambda^{3/2}}
\begin{pmatrix}
\mathrm{e}^{-\sqrt{\lambda}\,v}&0
\\
-\mathrm{e}^{\sqrt{\lambda}\,v}&0
\end{pmatrix}
+\frac{\mu}{\lambda^2}\begin{pmatrix}
\sqrt{\lambda}+\nu & \nu\mathrm{e}^{-2\sqrt{\lambda}\,v}
\\
-\nu\mathrm{e}^{2\sqrt{\lambda}\,v} & \sqrt{\lambda}-\nu
\end{pmatrix}\!\!
\begin{pmatrix}
\mathrm{e}^{-\sqrt{\lambda}\,u}&0
\\
-\mathrm{e}^{\sqrt{\lambda}\,u}&0
\end{pmatrix}
\\
&
=\frac{1}{\lambda^2}
\begin{pmatrix}
\nu\big(\sqrt\lambda\,\mathrm{e}^{-\sqrt{\lambda}\,v}
-\mu\mathrm{e}^{\sqrt{\lambda}\,(u-2v)}\big)
+\big(\sqrt\lambda+\nu\big)\mu\mathrm{e}^{-\sqrt{\lambda}\,u}
&0
\\
-\nu\big(\sqrt\lambda\,\mathrm{e}^{\sqrt{\lambda}\,v}
-\mu\mathrm{e}^{\sqrt{\lambda}\,(2v-u)}\big)
-\big(\sqrt\lambda-\nu\big)\mu\mathrm{e}^{\sqrt{\lambda}\,u}
&0
\end{pmatrix}\!.
\end{align*}
Hence
\begin{align*}
\bar{\gamma}_0
&
=\frac{\nu\big(\sqrt\lambda\,\mathrm{e}^{-\sqrt{\lambda}\,v}
-\mu\mathrm{e}^{\sqrt{\lambda}\,(u-2v)}\big)
+\big(\sqrt\lambda+\nu\big)\mu\mathrm{e}^{-\sqrt{\lambda}\,u}}
{\lambda\big[\big(\sqrt{\lambda}+\mu\big)\!\big(\sqrt{\lambda}+\nu\big)
-\mu\nu \mathrm{e}^{2\sqrt{\lambda}\,(u-v)}\big]}
\\
&
=\frac{\mu\big[\sqrt\lambda+\nu\big(1-\mathrm{e}^{\sqrt{\lambda}\,(u-v)}\big)\big]
\mathrm{e}^{-\sqrt{\lambda}\,u}
+\nu\big[\sqrt\lambda+\mu\big(1-\mathrm{e}^{\sqrt{\lambda}\,(u-v)}\big)\big]
\mathrm{e}^{-\sqrt{\lambda}\,v}}
{\lambda\big[\big(\sqrt{\lambda}+\mu\big)\!\big(\sqrt{\lambda}+\nu\big)
-\mu\nu \mathrm{e}^{2\sqrt{\lambda}\,(u-v)}\big]}.
\end{align*}
Therefore, since
$\bar{B}_0=-\bar{\gamma}_0\begin{pmatrix}1\\0\end{pmatrix}$,
we get
$$
\big(1\;\;0\big) N(x)\bar{B}_0=-\frac{\mu\big[\sqrt\lambda
+\nu\big(1-\mathrm{e}^{\sqrt{\lambda}\,(u-v)}\big)\big]
\mathrm{e}^{\sqrt{\lambda}\,(x-u)}
+\nu\big[\sqrt\lambda+\mu\big(1-\mathrm{e}^{\sqrt{\lambda}\,(u-v)}\big)\big]
\mathrm{e}^{\sqrt{\lambda}\,(x-v)}}
{\lambda\big[\big(\sqrt{\lambda}+\mu\big)\!\big(\sqrt{\lambda}+\nu\big)
-\mu\nu \mathrm{e}^{2\sqrt{\lambda}\,(u-v)}\big]}.
$$
Now,
$$
\bar{B}_1=\bar{S}_1\!\begin{pmatrix}1\\0\end{pmatrix}
-\bar{\gamma}_0\bar{R}_1\!\begin{pmatrix}1\\0\end{pmatrix}
=\frac{1}{\lambda^{3/2}}\begin{pmatrix}
\mu\mathrm{e}^{-\sqrt{\lambda}\,u}-\bar{\gamma}_0\big(\sqrt{\lambda}+\mu\big)
\\
-\mu\mathrm{e}^{\sqrt{\lambda}\,u}+\bar{\gamma}_0\mu\mathrm{e}^{2\sqrt{\lambda}\,u}
\end{pmatrix}
$$
Straightforward computations show that numerators of the entries of $\bar{B}_1$ can be
simplified into:
\begin{align*}
\mu\mathrm{e}^{-\sqrt{\lambda}\,u}+\bar{\gamma}_0\mu\mathrm{e}^{2\sqrt{\lambda}\,u}
&
=-\frac{\nu\sqrt\lambda\,\big[\sqrt\lambda+\mu\big(1-\mathrm{e}^{\sqrt{\lambda}\,(u-v)}\big)\big]
\mathrm{e}^{-\sqrt{\lambda}\,v}}
{\big(\sqrt{\lambda}+\mu\big)\!\big(\sqrt{\lambda}+\nu\big)
-\mu\nu \mathrm{e}^{2\sqrt{\lambda}\,(u-v)}},
\\
-\mu\mathrm{e}^{\sqrt{\lambda}\,u}+\bar{\gamma}_0\big(\sqrt{\lambda}-\mu\big)
&
=-\frac{\mu\sqrt\lambda\,\big[\sqrt\lambda+\nu\big(1-\mathrm{e}^{\sqrt{\lambda}\,(u-v)}\big)\big]
\mathrm{e}^{\sqrt{\lambda}\,u}}
{\big(\sqrt{\lambda}+\mu\big)\!\big(\sqrt{\lambda}+\nu\big)
-\mu\nu \mathrm{e}^{2\sqrt{\lambda}\,(u-v)}}.
\end{align*}
Thus
$$
\bar{B}_1=-\frac{1}{\lambda\big[\big(\sqrt{\lambda}+\mu\big)\!\big(\sqrt{\lambda}+\nu\big)
-\mu\nu \mathrm{e}^{2\sqrt{\lambda}\,(u-v)}\big]}
\begin{pmatrix}
\nu\big[\sqrt\lambda+\mu\big(1-\mathrm{e}^{\sqrt{\lambda}\,(u-v)}\big)\big]
\mathrm{e}^{-\sqrt{\lambda}\,v}
\\
\mu\big[\sqrt\lambda+\nu\big(1-\mathrm{e}^{\sqrt{\lambda}\,(u-v)}\big)\big]
\mathrm{e}^{\sqrt{\lambda}\,u}
\end{pmatrix}
$$
and then
$$
\big(1\;\;0\big) N(x)\bar{B}_1=
-\frac{\nu\big[\sqrt\lambda+\mu\big(1-\mathrm{e}^{\sqrt{\lambda}\,(u-v)}\big)\big]
\mathrm{e}^{\sqrt{\lambda}\,(x-v)}
+\mu\big[\sqrt\lambda+\nu\big(1-\mathrm{e}^{\sqrt{\lambda}\,(u-v)}\big)\big]
\mathrm{e}^{\sqrt{\lambda}\,(u-x)}
}{\lambda\big[\big(\sqrt{\lambda}+\mu\big)\!\big(\sqrt{\lambda}+\nu\big)
-\mu\nu \mathrm{e}^{2\sqrt{\lambda}\,(u-v)}\big]}.
$$
Furthermore,
\begin{align*}
\lqn{\bar{B}_2=\bar{S}_2\!\begin{pmatrix}1\\0\end{pmatrix}
-\bar{\gamma}_0\bar{R}_2\!\begin{pmatrix}1\\0\end{pmatrix}}
=\frac{1}{\lambda^{3/2}}\begin{pmatrix}
\nu\big(\sqrt\lambda\,\mathrm{e}^{-\sqrt{\lambda}\,v}
-\mu\mathrm{e}^{\sqrt{\lambda}\,(u-2v)}\big)
+\big(\sqrt\lambda+\nu\big)\mu\mathrm{e}^{-\sqrt{\lambda}\,u}
-\bar{\gamma}_0\big[\big(\sqrt{\lambda}+\mu\big)\!\big(\sqrt{\lambda}+\nu\big)
-\mu\nu \mathrm{e}^{2\sqrt{\lambda}\,(u-v)}\big]
\\
-\nu\big(\sqrt\lambda\,\mathrm{e}^{\sqrt{\lambda}\,v}
-\mu\mathrm{e}^{\sqrt{\lambda}\,(2v-u)}\big)
-\big(\sqrt\lambda-\nu\big)\mu\mathrm{e}^{\sqrt{\lambda}\,u}
+\bar{\gamma}_0\big[\big(\sqrt{\lambda}+\nu\big)\mu\mathrm{e}^{2\sqrt{\lambda}\,u}
-\big(\sqrt{\lambda}+\mu\big)\nu\mathrm{e}^{2\sqrt{\lambda}\,v}\big]
\end{pmatrix}\!.
\end{align*}
Obviously, the first entry of $\bar{B}_2$ vanishes.
Tedious computations show that the numerator of the second entry of $\bar{B}_2$ can be
simplified into:
\begin{align*}
-\mu\mathrm{e}^{\sqrt{\lambda}\,u}+\bar{\gamma}_0\big(\sqrt{\lambda}-\mu\big)
&
=-\frac{\mu\sqrt\lambda\,\big[\sqrt\lambda+\nu\big(1-\mathrm{e}^{\sqrt{\lambda}\,(u-v)}\big)\big]
\mathrm{e}^{\sqrt{\lambda}\,u}}
{\big(\sqrt{\lambda}+\mu\big)\!\big(\sqrt{\lambda}+\nu\big)
-\mu\nu \mathrm{e}^{2\sqrt{\lambda}\,(u-v)}}.
\end{align*}
Thus $\bar{B}_2=-\bar{\delta}_2\begin{pmatrix}0\\1\end{pmatrix}$
with
$$
\bar{\delta}_2=\frac{\mu\big[\sqrt\lambda
+\nu\big(1-\mathrm{e}^{\sqrt{\lambda}\,(u-v)}\big)\big]
\mathrm{e}^{\sqrt{\lambda}\,u}
+\nu\big[\sqrt\lambda+\mu\big(1-\mathrm{e}^{\sqrt{\lambda}\,(u-v)}\big)\big]
\mathrm{e}^{\sqrt{\lambda}\,v}
}{\lambda\big[\big(\sqrt{\lambda}+\mu\big)\!\big(\sqrt{\lambda}+\nu\big)
-\mu\nu \mathrm{e}^{2\sqrt{\lambda}\,(u-v)}\big]}
$$
and then
$$
\big(1\;\;0\big) N(x)\bar{B}_2=-\frac{\mu\big[\sqrt\lambda
+\nu\big(1-\mathrm{e}^{\sqrt{\lambda}\,(u-v)}\big)\big]
\mathrm{e}^{\sqrt{\lambda}\,(u-x)}
+\nu\big[\sqrt\lambda+\mu\big(1-\mathrm{e}^{\sqrt{\lambda}\,(u-v)}\big)\big]
\mathrm{e}^{\sqrt{\lambda}\,(v-x)}}
{\lambda\big[\big(\sqrt{\lambda}+\mu\big)\!\big(\sqrt{\lambda}+\nu\big)
-\mu\nu \mathrm{e}^{2\sqrt{\lambda}\,(u-v)}\big]}.
$$

As a result, we can see that the iterated Laplace transform of the couple of local times
$(L_t^u,L_t^v)$ at $u$ and $v$ can be expressed by the unified formula, for any $x\in \mathbb{R}$,
\begin{align*}
\lqn{\int_0^{\infty} \mathrm{e}^{-\lambda t}\, \mathbb{E}_x\!
\left(\mathrm{e}^{-\mu L_t^u-\nu L_t^v}\right) \mathrm{d} t}
=\frac{1}{\lambda}\left[1-
\frac{\mu\big[\sqrt\lambda+\nu\big(1-\mathrm{e}^{\sqrt{\lambda}\,(u-v)}\big)\big]
\mathrm{e}^{-\sqrt{\lambda}\,|x-u|}
+\nu\big[\sqrt\lambda+\mu\big(1-\mathrm{e}^{\sqrt{\lambda}\,(u-v)}\big)\big]
\mathrm{e}^{-\sqrt{\lambda}\,|x-v|}
}{\big(\sqrt{\lambda}+\mu\big)\!\big(\sqrt{\lambda}+\nu\big)
-\mu\nu \mathrm{e}^{2\sqrt{\lambda}\,(u-v)}}\right]\!.
\end{align*}
We retrieve formula~(1.18.1), p. 150 of~\cite{borodin}.
Consequently, the iterated Laplace transform of the local time $L_t^{u,v}=L_t^u+L_t^v$
in $\{u,v\}$ is given, for any $x\in \mathbb{R}$, by
$$
\int_0^{\infty} \mathrm{e}^{-\lambda t}\, \mathbb{E}_x\!
\left(\mathrm{e}^{-\mu L_t^{u,v}}\right) \mathrm{d} t
=\frac{1}{\lambda}\left[1-\mu\,
\frac{\big[\sqrt\lambda+\mu\big(1-\mathrm{e}^{\sqrt{\lambda}\,(u-v)}\big)\big]
\big[\mathrm{e}^{-\sqrt{\lambda}\,|x-u|}+\mathrm{e}^{-\sqrt{\lambda}\,|x-v|}\big]
}{\big(\sqrt{\lambda}+\mu\big)^2
-\mu^2 \mathrm{e}^{2\sqrt{\lambda}\,(u-v)}}\right]\!.
$$


\end{document}